\documentclass[11pt,a4paper,final]{article}
\pdfoutput=1
\usepackage{amssymb,latexsym}

\input{epsf.sty}
\usepackage{epsfig}
\usepackage{graphicx,color}
\usepackage{cite}
\usepackage{epstopdf}
\usepackage[top=1in, bottom=1in, left=1in, right=1in]{geometry}

\usepackage{graphics}
\usepackage{subfigure}
\usepackage{setspace}

\usepackage{amsmath}

\usepackage{mathrsfs,amsfonts}

\usepackage{amsthm,amsxtra}

\usepackage[final]{showkeys}

\usepackage{stmaryrd}

\usepackage{algorithm}
\usepackage{algorithmicx}
\usepackage{algpseudocode}
\usepackage{mathtools}
\newif\ifPDF
\ifx\pdfoutput\undefined
\PDFfalse
\else
\ifnum\pdfoutput > 0
\PDFtrue
\else
\PDFfalse
\fi
\fi

\ifPDF
\usepackage{pdftricks}
\begin{psinputs}
	\usepackage{pstricks}
	\usepackage{pstcol}
	\usepackage{pst-plot}
	\usepackage{pst-tree}
	\usepackage{pst-eps}
	\usepackage{multido}
	\usepackage{pst-node}
	\usepackage{pst-eps}
\end{psinputs}
\else
\usepackage{pstricks}
\fi

\ifPDF
\usepackage[debug,pdftex,colorlinks=true,
linkcolor=blue, bookmarksopen=false,
plainpages=false,pdfpagelabels]{hyperref}
\else
\usepackage[dvips]{hyperref}
\fi

\usepackage[openbib]{currvita}

\usepackage{fancyhdr}

\newtheorem{theorem}{Theorem}[section]

\newtheorem{remark}[theorem]{Remark}

\renewcommand\appendix{\par
  \setcounter{section}{0}
  \setcounter{subsection}{0}
  \setcounter{figure}{0}
  \setcounter{table}{0}
  \renewcommand\thesection{Appendix \Alph{section}}
  \renewcommand\theequation{\Alph{section}.\arabic{equation}}
  \renewcommand\thefigure{\Alph{section}.\arabic{figure}}
  \renewcommand\thetable{\Alph{section}.\arabic{table}}
  \renewcommand\thethm{\Alph{section}.\arabic{thm}}
}

\numberwithin{equation}{section}

\date{}

\title{A finite element contour integral method for computing the resonances of metallic grating structures with subwavelength holes}

\author{Yingxia Xi \thanks{School of Mathematics and Statistics, Nanjing University of Science and Technology, Nanjing, 210094, China. \tt xiyingxia@njust.edu.cn.}
\and Junshan Lin \thanks{\footnotesize Department of Mathematics and Statistics, Auburn University, Auburn, AL 36849. \tt jzl0097@auburn.edu.}
\and Jiguang Sun \thanks{Department of Mathematical Sciences, Michigan Technological University, Houghton, MI 49931. \tt jiguangs@mtu.edu.}
}

\begin{document}
\maketitle

\begin{abstract}
We consider the numerical computation of resonances for metallic grating structures with dispersive media and small slit holes. The underlying eigenvalue problem is nonlinear and the mathematical model is multiscale due to the existence of several length scales in problem geometry and material contrast.
We discretize the partial differential equation model over the truncated domain using the finite element method and develop a multi-step contour integral eigensolver to compute the resonances. The eigensolver first locates eigenvalues using a spectral indicator and then computes eigenvalues by a subspace projection scheme. The proposed numerical method is robust and scalable, and does not require initial guess as the iteration methods.
Numerical examples are presented to demonstrate its effectiveness.
\end{abstract}

\setcounter{equation}{0}
\setlength{\arraycolsep}{0.25em}

\section{Introduction}
Resonances play a significant role in the design of novel materials, due to their ability to generate unusual physical phenomena that open up a broad possibility in modern science and technology. Typically the resonances could be induced by arranging the material parameters or the structure geometry carefully with high-precision fabrication techniques available nowadays.
Mathematically, resonances correspond to certain complex eigenvalues of the underlying differential operators with the corresponding eigenmodes that are either localized with finite energy or extended to infinity.
When the resonances are excited by external wave field at the resonance frequencies, the wave field generated by the system can be significantly amplified, which leads
to various important applications in acoustics and electrodynamics, etc.

One important class of resonant optical materials is the subwavelength nano-holes perforated in noble metals, such as gold or silver. Tremendous research has been sparked in the past two decades in pursuit of more efficient resonant nano-hole devices (cf. \cite{garcia10, rodrigo16} and references therein) since the seminal work  \cite{ebbesen98}. At the resonant frequencies, the optical transmission through the tiny holes exhibit extraordinary large values, or the so-called extraordinary optical transmission (EOT), which can be used for biological and chemical sensing, and the design of novel optical devices, etc \cite{blanchard17, cetin2015, huang08, li17plasmonic, rodrigo16}.
The main mechanisms for the EOT in the subwavelength hole devices are resonances. These include scattering resonances induced by the tiny holes patterned in the structure and surface plasmonic resonances generated from the metallic materials \cite{garcia10}. Although both
are eigenvalues of the differential operator when formulated in a finite domain, their eigenmodes are very different.


While significant progress has been made on the mathematical studies of resonances in such subwavelength metallic structures, the studies are mostly based on the ideal models when the metal is a perfect conductor and the shape of the hole is simple \cite{bontri10, bontri10-2, fatima21, holsch19,  liazou20, linluzha23, linshizha20, linzha17, linzha18a, linzha18b, linzha21, zhlu21,  luwanzho21}. Such an assumption allows one to impose the boundary conditions over the boundary of the metal in the mathematical model and analyze the resonances induced by the subwavelength holes patterned in the metal. However, the model neglects the penetration of the wave field into the metal, which is significant in most optical and acoustic frequency regime \cite{maier07}. Hence the second type of resonances, namely the surface plasmonic resonance which is significant in metallic structures, is absent in the ideal model.

In this paper, we consider the more challenging model for which the permittivity of the metal is described by a frequency-dependent function and the shape of the hole can be arbitrary.
We develop a finite element contour integral approach to compute the resonances of the multiscale metallic structure \cite{SunZhou2016, Beyn2012, Huang2016, HuangSunYang, XiaoSun2021JOSAA, Gong2022MC}.
Since wave can penetrate into the metal, we consider the full transmission problem for the Helmholtz equation in which the permittivity function of the material is defined piecewisely and depends on frequency.
While there exist many works for computing the related scattering problems (see, for instance \cite{astilean, kriegsmann, porto}), the research on solving the corresponding eigenvalue problems is scarce.

Two major challenges arise when solving the eigenvalue problem numerically:
\begin{itemize}

\item[(i)] There exist several length scales in the problem geometry and the material contrast between the metal and background could be large. Typically
the size of the tiny hole and the skin depth characterizing the wave penetration depth into the metal are much smaller than the free-space wavelength. This requires resolving the wave oscillation
at fine scales accurately. In addition, the permittivity contrast between metal and the background medium could exceed 100 in certain frequency regime. We truncate the problem into a finite domain and employ a finite element discretization for the differential operator with unstructured meshes to resolve the wave oscillation accurately \cite{SunZhou2016}.

\item[(ii)] The dielectric function of the metal depends nonlinearly on the wave frequencies, or the eigen-parameters, as such the eigenvalue problem is nonlinear. In general, nonlinear eigenvalue solvers such as the Newton type methods require sufficiently close initial guesses to ensure the convergence, which are usually unavailable.
We design a robust contour integral method to locate and compute the eigenvalues of the discretized system. First, a spectral indicator method 
is used to locate the eigenvalues by examining the regions on the complex plane. Then the subspace projection scheme (cf. \cite{Beyn2012}) is employed to compute the eigenvalues accurately. Both methods rely on the contour integral of the resolvent operator. Finally, the verification of eigenvalues is performed by using the discretized algebraic system.
The proposed method is scalable since different regions on the complex plane can be examined in parallel.


\end{itemize}
The proposed computational framework is more versatile than the mode matching method developed in \cite{Lin} for solving the eigenvalue problem, which relies on the expansion of the wave inside the tiny whole and requires special shape of the hole geometry. It is also more flexible than the integral equation method in \cite{linzha19} as the evaluation of the Green's function in grating can be slow and complicated. Very importantly, the proposed method locate the eigenvalues using the spectral indicator and
does not need good initial guesses as required by Newton type methods.
We would like to point out that the considered eigenvalue problem is closely related to nanoparticle plasmonic resonance problem, in which the permittivity of the metal is also a frequency-dependent function but the problem is imposed over the finite-region nanoparticles \cite{ammari1, ammari2, ammari3, ammari4}. The configuration of the periodic structure is similar to dielectric grating or periodic crystal slab in general, for which the mathematical study has been restricted to dielectric materials \cite{bao, shipman}.

The rest of the paper is organized as follows. In Section~\ref{Problem setup}, we introduce the mathematical model and the Dirichelt-to-Neumann (DtN) map to truncate the infinite domain to a finite one. Section~\ref{FEM} presents the discrete weak formulation and the finite element scheme, which leads to a nonlinear algebraic eigenvalue problem. In Section~\ref{MCI}, we propose a multi-step eigensolver for the discrete system based on complex contour integrals and discuss the implementation details.
Numerical examples are presented in Section~\ref{NE} to test the effectiveness of the proposed method. We consider mathematical models with different permittivity functions and slit geometries. In particular, two types of resonances for the metallic structure are examined by the corresponding eigenfunctions.
Finally, we summarize the work and discuss future direction along this line in Section~\ref{CFW}.

\section{The eigenvalue problem}\label{Problem setup}
\begin{figure}[!htbp]
\begin{center}
\includegraphics[width=15cm]{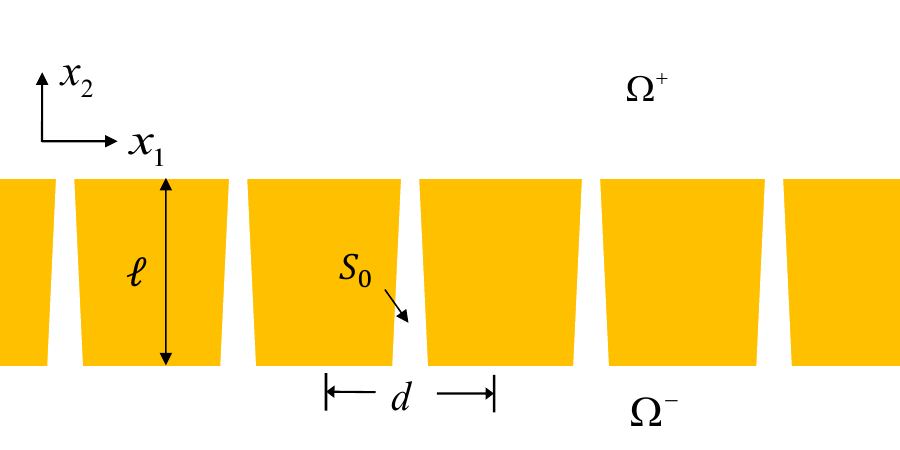}
\caption{Geometry of the periodic metallic structure.
The slits $S$ are arranged periodically with the size of the period $d$.
The domains above and below the metallic slab are denoted as $\Omega^{+}$ and $\Omega^{-}$ respectively, and the domain of the metal is denoted as $\Omega$. }\label{fig-prob_geo}
\end{center}
\end{figure}

We consider a metallic slab that is perforated with a periodic array of slits
and the geometry of its cross section is depicted in Figure \ref{fig-prob_geo}. The  slab occupies the domain
$\Omega_0:=\{(x_1,x_2)\;|\; -\frac{\ell}{2}<x_2<\frac{\ell}{2}\}$, and
the slits occupy the region
$\displaystyle{S=\bigcup_{n=0}^{\infty} (S_0 + nd)}$,
where $d$ is the size of the period and
$S_{0}$ represents the slit in the reference period.
Denote the domain of the metallic structure by $\Omega:=\Omega_0\backslash S$.
The semi-infinite domain above and below the slab is denoted by $\Omega^{+}$ and $\Omega^{-}$.
The relative electric permittivity $\varepsilon(\omega; x)$ on the $x_1x_2$ plane is given by
\begin{equation*}
\varepsilon(\omega; x)= \left\{
\begin{array}{lll}
\medskip
\varepsilon_0,  & x\in\Omega^+\cup\Omega^-, \\
\medskip
\varepsilon_m(\omega),  &x\in\Omega,
\end{array}
\right.
\end{equation*}
where $\varepsilon_0$ and $\varepsilon_m$ denote the relative permittivity in the vacuum and metal, respectively.  $\omega$ represents the operating frequency of the wave. The permittivity in the metal is frequency-dependent. In this work, we consider the so-called Drude model ( cf. \cite{ordal83} ) such that
\begin{equation}\label{eq:em}
\varepsilon_m(\omega)=1-\frac{\omega_p^2}{\omega^2+i\Gamma \omega},
\end{equation}
where $\omega_p$ is the volume plasma frequency and $\Gamma$ is the damping coefficient. For convenience of notation, we define the wavenumber $k=\omega/c$, wherein $c$ is the wave speed, and rewrite the permittivity function as
\begin{equation}
\varepsilon_m(k)=1-\frac{\bar\omega_p^2}{k^2+i\bar\Gamma k},
\end{equation}
where $\bar{\omega}_p=\frac{\omega_p}{c}$ and $\bar{\Gamma}=\frac{\Gamma}{c}$.

We consider the following eigenvalue problem
for the transverse magnetic (TM) polarized electromagnetic wave:
\begin{equation}\label{eq:Helmholtz}
    \nabla \cdot \left( \dfrac{1}{\varepsilon(k;x)} \nabla u \right) + k^2 u = 0   \quad\quad  \mbox{in} \; \mathbb{R}^2,
\end{equation}
where $u$ represents the $x_3$-component of the magnetic field.
In addition, along the metal boundary $\partial\Omega$, there holds
\begin{equation}\label{eq:interface_cond}
    [u] = 0,  \quad \left[\dfrac{1}{\varepsilon} \dfrac{\partial u}{\partial \nu}\right] = 0,
\end{equation}
wherein $[\cdot]$ denotes the jump of the quantity when the limit is taken along the positive and negative unit normal direction $\nu$ of $\partial\Omega$, respectively.

We restrict the problem to one periodic cell $\{ (x_1, x_2) \; | \; 0< x_1 <d  \}$ following the Floquet-Bloch theory \cite{kuchment1993}.
For each Bloch wavenumber $\kappa$ in the Brillouin zone $[-\pi/d, \pi/d]$, we look for quasi-periodic solutions
of \eqref{eq:Helmholtz}-\eqref{eq:interface_cond} such that $u(x_1,x_2)=e^{i\kappa x_1}\tilde u(x_1,x_2)$,
where $\tilde u$ is a periodic function with $\tilde u(x_1+d,x_2)=\tilde u(x_1,x_2)$.
This gives the so-called quasi-periodic boundary condition for $u$ over the periodic cell
\begin{equation}\label{eq-quasi_periodic}
u(d,x_2)=e^{i\kappa d} u(0,x_2), \quad \partial_{x_1}u(d,x_2)=e^{i\kappa d} \partial_{x_1} u(0,x_2).
\end{equation}
Along the $x_2$ direction, we assume that the wave is outgoing. Then using the Fourier expansion, it can be shown that
$u$ adopts the following Rayleigh-Bloch expansion
\begin{equation}\label{eq-rad_cond}
 u(x_1,x_2) =  \sum_{n=-\infty}^{\infty} u_n^{+} e^{i \kappa_n x_1 + i\zeta_n  (x_2-H) } \quad \mbox{and} \quad
 u(x_1,x_2) =  \sum_{n=-\infty}^{\infty} u_n^{-} e^{i \kappa_n x_1 - i\zeta_n  (x_2+H) }
\end{equation}
in the domain $\{x_2>H\}$ and $\{x_2<-H\}$, respectively. Here $H$ is some positive constant satisfying $H>\frac{\ell}{2}$,
$$
\kappa_n=\kappa+\dfrac{2\pi n}{d} \quad \mbox{and} \quad
\zeta_n(k)= \sqrt{k^2-\kappa_n^2},
$$
are the wavenumbers in the $x_1$ and $x_2$ directions,
and the function $f(z)=\sqrt{z}$ is understood as an analytic function defined in the domain $\mathbf{C}\backslash\{-it: t\geq 0\}$ by
$$
z^{\frac{1}{2}} = |z|^{\frac{1}{2}} e^{\frac{1}{2} i\arg z}.
$$
The coefficients $u_n^{\pm}$ are the Fourier coefficients of the solution at $x_2=\pm H$ such that
$$ u(x_1,\pm H) =  \sum_{n=-\infty}^{\infty} u_n^{\pm} e^{i \kappa_n x_1}.  $$
The rest of paper is devoted to solving for the eigenvalue $k$ for the problem \eqref{eq:Helmholtz}-\eqref{eq-rad_cond}.

To reduce the problem to a bounded domain, one can take the normal derivative of the Rayleigh-Bloch expansion above and define the Dirichlet-to-Neumann map at $x_2=\pm H$:
\begin{eqnarray*}
&& T^+\Big[u(\cdot,H)\Big] :=  \sum_{n=-\infty}^{\infty} i\zeta_n u_n^{+} e^{i \kappa_n x_1},\\
&& T^-\Big[u(\cdot,-H)\Big]  :=  -\sum_{n=-\infty}^{\infty} i\zeta_n u_n^{-} e^{i \kappa_n x_1}.
\end{eqnarray*}
Then for each $\kappa \in [-\pi/d, \pi/d]$, the eigenvalue problem can be formulated as a nonlinear eigenvalue problem in the bounded domain $\Omega_H:=\{ (x_1,x_2) \,|\, 0<x_1<d, -H<x_2<H \}$ as follows:
\begin{equation} \label{eq:eig_bdd_domain}
\left\{
\begin{array}{llll}
\vspace*{0.3cm}
\nabla \cdot \left( \dfrac{1}{\varepsilon(k; x)} \nabla u \right) + k^2 u = 0   \quad\quad  \mbox{in} \; \Omega_H,  \\
\vspace*{0.1cm}
u(d,x_2)=e^{i\kappa d} u(0,x_2), \quad \partial_{x_1}u(d,x_2)=e^{i\kappa d} \partial_{x_1} u(0,x_2),    \\
\vspace*{0.1cm}
\dfrac{\partial u}{\partial \nu}(x_1,\pm H) = T^\pm \Big[u(\cdot,\pm H)\Big]   \quad\quad  \mbox{on} \; x_2 = \pm H.
\end{array}
\right.
\end{equation}
The correponding weak formulation is to find $k \in \mathbb C$ and $u \in  H_\kappa^1(\Omega_H)$ such that
\begin{equation}\label{eq:weak_main}
\int_{\Omega_H} \dfrac{1}{\varepsilon(k, x)} \nabla u \nabla \bar v - k^2 u \bar v dx - \langle  T^+ u, v \rangle - \langle  T^- u, v \rangle = 0 \quad \text{for all }v\in H_\kappa^1(\Omega_H),
\end{equation}
where $H_\kappa^1(\Omega_H)$ is the Sobolev space defined by $$H_\kappa^1(\Omega_H) = \{ u \in L^2(\Omega_H) | \partial_{x_j} u  \in L^2(\Omega_H), u(d,x_2)=e^{i\kappa d} u(0,x_2), \partial_{x_1}u(d,x_2)=e^{i\kappa d} \partial_{x_1} u(0,x_2  \}. $$

\begin{remark}
Due to the Rayleigh-Bloch expansion \eqref{eq-rad_cond}, the eigenvalue problem \eqref{eq:eig_bdd_domain} is nonlinear even if $\varepsilon$ does not depend on $k$.
\end{remark}


\section{Finite element discretization}\label{FEM}
We employ a finite element method to discretize the weak formulation \eqref{eq:weak_main}.
Let $\Omega_H$ be covered with a regular and quasi-uniform mesh $\mathcal{T}_h$ consisting of triangular elements $\{T_j\}_{j=1}^{N_e}$. The mesh size is defined as $h=\mathop{max}\limits_{1\leq j\leq N_e} \rho_j$, where $\rho_j$ is the diameter of the inscribed circle of $T_j$.
Denote by $V_h$ the linear Lagrange finite element space associated with $\mathcal{T}_h$. The subspace $V_{\kappa,h} \subset V_h$ contains the functions that satisfy the quasi-periodic boundary condition \eqref{eq-quasi_periodic}. In addition, we denote by $V_h^{B\text{-}upper}$ the subspace of $V_h$ with degrees of freedom (DOF) on $x_2 = H$ and $V_h^{B\text{-}lower}$ the subspace of $V_h$ with degrees of freedom (DOF) on $x_2 = -H$.





On the discrete level, we truncate the infinite series of the DtN mappings $T^\pm u$:
\begin{eqnarray*}
&&T_{D_t}^+\left(u^s(x_1,H)\right) =: \sum_{n=-D_t}^{D_t} i\zeta_n u_n^{s,+} e^{i \kappa_n x_1}, \\
&&T_{D_t}^-\left(u^s(x_1,-H)\right) =: \sum_{n=-D_t}^{D_t} i\zeta_n u_n^{s,-} e^{i \kappa_n x_1}.
\end{eqnarray*}
The non-negative integer $D_t$ is called the truncation order of the DtN mapping. The discrete problem for \eqref{eq:weak_main} is to find $k\in\mathbb{C}$ and $u_{h}\in V_{\kappa,h}$ such that
\begin{equation}\label{discrete_weak_formulation}
\int_{\mathcal{T}_h} \dfrac{1}{\varepsilon(k, x)} \nabla u_{h} \nabla \bar v_{h} - k^2 u_{h} \bar v_h dx - \langle  T_{D_t}^+ u_{h}, v_h \rangle - \langle  T_{D_t}^- u_{h}, v_h \rangle =0 \quad v_h\in V_{\kappa,h}.
\end{equation}

We now derive the matrix form for \eqref{discrete_weak_formulation}.
Assume that the basis funcitions for $V_h^{B\text{-}upper}$  and $V_h^{B\text{-}lower}$ are given by $\phi_1,\phi_2,\cdots,\phi_{N_b}$ and $\phi_{N_b+1},\phi_{N_b+2},\cdots,\phi_{2N_b}$, respectively,
and the basis functions for $V_h$ are given by $\phi_1,\cdots,\phi_{2N_b},\phi_{2N_b+1},\cdots,\phi_N$. 
Writing $u_h\in V_{\kappa,h}\subseteq V_h$ as
$
u_h=\sum_{j=1}^{N} u_j\phi_j,
$
the stiffness matrix $(A_1)_{N\times N}$ and mass matrix $(A_2)_{N\times N}$ are given by
$$(A_1)_{q,j}=\int\frac{1}{\varepsilon(k, x)}\nabla \phi_j\nabla \bar{\phi}_q dx,\ \ \
(A_2)_{q,j}=\int \phi_j \bar{\phi}_q dx,\ \ \ q,j=1,2,\cdots,N.
$$

Since the permittivity function $\varepsilon(k, x)$ depends nonlinearly on the wavenumber $k$, $A_1$ is a nonlinear matrix function of $k$, which is denoted by $A_1(k)$ for convenience.
Similarly, for a fixed wavenumber $k\in\mathbb{C}$, the matrices for $\langle  T_{D_t}^+ u_{h}, v_h \rangle$ and $\langle  T_{D_t}^- u_{h}, v_h \rangle$ are denoted as $(A_3(k))_{N\times N}$ and $(A_4(k))_{N\times N}$, respectively. The non-zero elements are given by
\[
(A_3(k))_{q,j}=\int_{\partial\Omega_H^+}\sum_{n=-D_t}^{D_t}i\zeta_n \left(\frac{1}{d}\int_{-\frac{d}{2}}^{\frac{d}{2}}\phi_je^{-i\kappa_n x_1}dx_1 \right)e^{i\kappa_n x_1}\bar{\phi}_q ds,
\]
where $q,j=1,2,\cdots,N_b$, and
\[
(A_4(k))_{q,j}=\int_{\partial\Omega_H^-}\sum_{n=-D_t}^{D_t}i\zeta_n \left(\frac{1}{d}\int_{-\frac{d}{2}}^{\frac{d}{2}}\phi_je^{-i\kappa_n x_1}dx_1 \right)e^{i\kappa_n x_1}  \bar{\phi}_q ds,
\]
where $q,j=N_b+1,N_b+2,\cdots,2N_b$.


The quasi-periodic boundary conditions (\ref{eq-quasi_periodic}) are treated using the Lagrange multiplier.
We illustrate how to enforce the condition
\begin{equation}\label{quaisi-periodic1}
u(d,x_2)=e^{i\kappa d} u(0,x_2).
\end{equation}
Note that the quasi-periodic boundary conditions for the partial derivative are satisfied naturally \cite{SukumarPask}.
Assume a one-to-one correspondence between the mesh nodes of $\mathcal{T}_h$ on the left boundary and the right boundary of $\Omega_H$. Let the nodes on the left boundary be $x^l_1,x^l_2,\cdots,x^l_J$. The corresponding nodes on the right boundary are $x^r_1, x^r_2,\cdots,x^r_J$. The condition (\ref{quaisi-periodic1}) implies that
\begin{equation}
e^{i\kappa d}u(x^l_1)=u(x^r_1),\ e^{i\kappa d}u(x^l_2)=u(x^r_2),\cdots,e^{i\kappa d}u(x^l_J)=u(x^r_J).
\end{equation}
Introduce the Lagrange multiplier $\hat{\lambda}=(\lambda_1,\lambda_2,\cdots,\lambda_J)^T$ and let $u=(u_1,u_2,\cdots,u_N)^T$. Define the auxiliary matrix $B_{J\times N}$ with non-zero elements given by
\begin{align} \nonumber
B(1,r_1)&=1,\ B(1,l_1)=-e^{i\kappa d},\\ \nonumber
B(2,r_2)&=1,\ B(2,l_2)=-e^{i\kappa d},\\ \nonumber
&\cdots\\ \nonumber
B(J,r_J)&=1,\ B(J,l_J)=-e^{i\kappa d},\nonumber
\end{align}
where $l_1, \ldots, l_J$ are the indices of the nodes $x^l_1, \ldots, x^l_J$ and $r_1, \ldots, r_J$ are the indices of the nodes $x^r_1, \ldots, x^r_J$.

The augmented algebraic eigenvalue system for \eqref{discrete_weak_formulation} is
\begin{equation}\label{Nonlinear_matrix}
\mathbb{G}(k)\xi\triangleq\left(
\begin{array}{cc}
A(k) &B^T \\
B &O_{J\times J}
\end{array}
\right)
\left(
\begin{array}{c}
 u \\
\hat{\lambda}
\end{array}
\right)= O_{(N+J)\times 1},
\end{equation}
where
\[
A(k)=A_1(k)-k^2A_2-A_3(k)-A_4(k)
\]
and $O_{J\times J}$, $O_{(N+J)\times 1}$ are null matrices. 
$\mathbb{G}(\cdot):\Omega\rightarrow\mathbb{C}^{(N+J),(N+J)}$ is a nonlinear matrix-valued function.
Let $v=(u;\hat{\lambda})$.
We call $k \in \mathbb{C}$ an eigenvalue and $v \in \mathbb{C}^{(N+J)}$ the associated eigenvector of $\mathbb{G}(\cdot)$ if
\begin{equation}\label{Gkv}
\mathbb{G}(k) v = 0.
\end{equation}
The eigenvalues $k$ of $\mathbb{G}(\cdot)$ are resonances that we are looking for.

For the examples considered in this paper, the multiple length scales are treated using unstructured meshes with finer grids near the metal and slits and coarser grids for the homogeneous background (see Figure~\ref{pec_example}). For more complicated problems, one should incorporate more sophisticated basis functions and implement a global numerical formulation that couples these multiscale basis functions \cite{EfendievHou2009}.

\section{Multi-step eigensolver based on the contour integral}\label{MCI}
We propose a multi-step scheme to compute the eigenvalues of $\mathbb{G}(\cdot)$ inside a given bounded region on the complex plane. It consists of three steps: (1) detection using the spectral indicator method \cite{HuangSunYang}, (2) computation using the projection method \cite{Beyn2012}, and (3) verification.

The main ingredient of Steps (1) and (2) is the contour integral. Let $R \subset \mathbb{C}$ be a simply-connected bounded domain with piecewise smooth boundary. We call $R$ the region of interest and the goal is to compute all eigenvalues inside it. Assume that $\mathbb{G}(k)$ is holomorphic on $R$ and $\mathbb{G}(k)^{-1}$ exists for all $k \in \partial R$.
Define a projection operator $P: \mathbb{C}^{(N+J)} \to \mathbb{C}^{(N+J)}$  by
\begin{equation}\label{Pv}
Pv:=\dfrac{1}{2\pi i}\int_{\partial R} \mathbb{G}(k)^{-1}v \, dk, \quad v \in \mathbb{C}^{(N+J)}.
\end{equation}

In Step (1), we cover $R$ with small disks $\{R_i\}_{i=1}^I$ (see Figure \ref{fig:indicator_step1}) and use the above operator
to determine if $R_i$ contains eigenvalues.
If $\mathbb{G}(k)$ has no eigenvalues in $R_i$, then $\mathbb{G}(k)^{-1}$ is holomorphic on $\overline{R}_i$. By Cauchy's Theorem, $Pv = 0$ for any $v \in \mathbb{C}^{(N+J)}$, where the contour integral \eqref{Pv} is evaluated over $\partial R_i$.
On the other hand, if $\mathbb{G}(k)$ attains eigenvalues in $R_i$, $Pv \ne 0$ almost surely for a random vector $v \in \mathbb{C}^{(N+J)}$.
Computationally we choose a random vector $v$ and compute $P_i v$, $i=1, \ldots, I$. We use $\|P_i v\|$ as the indicator for $R_i$. If $\|P_i v\| = O(1)$, then $R_i$ contains eigenvalues. If $\|P_i v\| = o(1)$, there is no eigenvalues in $R_i$, which is then discarded.
\begin{remark}
Ideally each $R_i$ is small enough and contains a few eigenvalues. This is usually done by trial and error. In the implementation, we set the threshold value as $0.2$, i.e., if $\|P_i v\| \ge 0.2$, we save $R_i$ for Step (2). We refer the readers to \cite{HuangSunYang} for some discussions on this choice.
\end{remark}

In Step (2), given a (small) disk $R_i$, we use the subspace projection method in \cite{Beyn2012} to compute candidate eigenvalues inside the disk.
For convenience of notation, we denote $R_i$ by $R$ in the following discussions.
We present the algorithm for the case of simple eigenvalues.
Eigenvalues with multiplicity more than one can be treated similarly (Theorem 3.3, \cite{Beyn2012}).

Assume that there exist $M$ eigenvalues $\{k_j\}_{j=1}^M$ inside $R$ and no eigenvalues lie on $\partial R$.
Let $f:\partial R\rightarrow \mathbb{C}$ be any holomorphic function. 
Then one has that
\begin{align}\label{expansionT}
\frac{1}{2\pi i}\int_{\partial R} f(z)\mathbb{G}(z)^{-1}dz=\sum_{j=1}^{M}f(k_j)v_jw_j^H,
\end{align}
where $v_j, w_j$ are left and right eigenfunctions corresponding to $k_j$ such that $w_j^H\mathbb{G}'(k_j)v_j=1$.
Let $V=[v_1,v_2,\cdots,v_M]\in\mathbb{C}^{(N+J),M}$ and $W=[w_1,w_2,\cdots,w_M]\in\mathbb{C}^{(N+J),M}$. Then
\begin{equation}
\frac{1}{2\pi i}\int_{\partial R} f(z)\mathbb{G}(z)^{-1}dz=V\Sigma_fW^H,
\end{equation}
where $\Sigma_f=diag(f(k_1),\cdots,f(k_M))$. The following theorem explains how to compute the eigenvalues $k_j$'s.
\begin{theorem}\label{beyn_theorem}
(Theorem 3.1, \cite{Beyn2012})
Let $\tilde{V}\in\mathcal{C}^{(N+J),L_1}(M\leq L_1\leq N+J)$ be chosen randomly. In a generic sense, the volumn vectors of $\tilde{V}$ are linearly independent. Then, it holds that $\rm{rank}(W^H\tilde{V})=M$ and $\rm{rank}(V)=M$.
Assume that
\begin{align}\label{integrals1}
C_0&=\frac{1}{2\pi i}\int_{\partial R}\mathbb{G}(z)^{-1}\tilde{V}dz \in \mathbb{C}^{(N+J)\times L_1}\\ \label{integrals2}
C_1&=\frac{1}{2\pi i}\int_{\partial R}z\mathbb{G}(z)^{-1}\tilde{V}dz \in \mathbb{C}^{(N+J)\times L_1}
\end{align}
and the singular value decomposition
$$C_0=V_0\Sigma_0W_0^H,$$
where $V_0\in \mathcal{C}^{(N+J)\times M}$, $\Sigma_0=diag(\sigma_1,\sigma_2,\cdots,\sigma_M)$, $W_0\in \mathcal{C}^{L_1\times M}$. Then, the matrix
\begin{equation}\label{D}
D:=V_0^HC_1W_0\Sigma_0^{-1}\in\mathcal{C}^{M,M}
\end{equation}
is diagonalizable with eigenvalues $k_1,k_2,\cdots,k_M$.
\end{theorem}

As a consequence of the above theorem, one computes \eqref{integrals1} and \eqref{integrals2} using the trapezoidal rule to obtain $C_0$ and $C_1$, performs the singular value decomposition for $C_0$, and then calculate eigenvalues of $D$.
Then the eigenvalues of $\mathbb{G}(k)$ in $R$ coincide with the eigenvalues of $D$.
Write the parameterization for $\partial R$ as
$$
\psi(\theta)=z_0+re^{i\theta},\ \ \ \theta\in(0,2\pi],
$$
where $z_0$ is the center and $r$ is the radius. Taking the equidistant nodes $\theta_j=\frac{2\pi j}{N_t}$, $j=1,2,\cdots,N_t$, and using the trapezoid rule, we obtain the following approximations for (\ref{integrals1}) and (\ref{integrals2}), respectively,
\begin{align}\label{integrals11}
C_0\approx C_{0,N_t}&=\frac{1}{iN_t}\sum_{j=1}^{N_t}\mathbb{G}(\psi(\theta_j))^{-1}\tilde{V}\psi'(\theta_j), 
\\ \label{integrals21} 
C_1\approx C_{1,N_t}&=\frac{1}{iN_t}\sum_{j=1}^{N_t}\mathbb{G}(\psi(\theta_j))^{-1}\tilde{V}\psi(\theta_j)\psi'(\theta_j). 
\end{align}

The number of eigenvalues inside $R$, i.e., $M$, is not known a priori. One would expect that there would be a gap between the group of large singular values of  $C_0$ and the group of small singular values of $C_0$.
However, this is not the case for the challenging problems considered in this paper. From the numerical examples, we observe that there is no significant gap between the singular values and one has to decide how many singular values to keep.  For robustness, we set a small tolerance value $\sigma_0$ such that if there are $M_0$ singular values of $C_0$ that are larger than $\sigma_0$, we compute $M_0$ eigenvalues of $D$ as the output values of Step (2).

In Step (3), we substitute the output values from Step (2) into \eqref{Gkv} to obtain $M_0$ matrices $\mathbb{G}(k_i), i= 1, \ldots, M_0,$ and compute the smallest eigenvalue $\lambda_{k_i}^0$ of $\mathbb{G}(k_i)$. If $\lambda_{k_i}^0 < 10^{-12}$, $k_i$ is taken as an eigenvalue of $\mathbb{G}(\cdot)$. If $\lambda_{k_i}^0 > 10^{-5}$, $k_i$ is discarded. If $\lambda_{k_i}^0$ is such that $10^{-12} \le \lambda_{k_i}^0 \le 10^{-5}$, an additional round of computation is performed. The values $10^{-12}$ and $10^{-5}$ are problem dependent and chosen by trial and error. Typical eigensovlers such as Arnodi methods can be used to compute the smallest eigenvalue of $\mathbb{G}(k_i)$.

The following algorithm summarizes the multi-step contour integral method to compute the eigenvalues for $\mathbb{G}(\cdot)$ in $R$.

\vskip 0.2cm
\textbf{Algorithm 1}
\begin{itemize}
\item[-] Given a region $R\subset \mathbb C$ of interest, compute eigenvalues of $\mathbb{G}(\cdot)$ in $R$.
\item[(1)] Identify small sub-regions of $R$ that might contains eigenvalues
    \begin{itemize}
    \item[(1.a)] Cover $R$ by smaller disks $R_i, i=1, \ldots I,$ and pick a random vector $v$.
    \item[(1.b)] Compute and normalize $\|P_iv\|$. Store $R_i$'s such that $\|P_iv\| \ge 0.2$.
    \end{itemize}
\item[(2)] For each stored $R_i$, compute the candidate eigenvalues.
    \begin{itemize}
    \item[(2.a)] Choose a large enough $L_1\leq N+J$ and generate a random matrix $\tilde{V}\in \mathcal{C}^{(N+J),L_1}$.
    \item[(2.b)] Calculate numerical integration \eqref{integrals11} and \eqref{integrals21}.
    \item[(2.c)]  Compute the singular value decomposition
            $$C_{0,N_t}=V\Sigma W^H,$$
        where $V\in\mathbb{C}^{(N+J),L_1}$, $\Sigma=diag(\sigma_1,\sigma_2,\cdots,\sigma_{L_1})$, $W\in\mathbb{C}^{L_1,L_1}$.
    \item[(2.d)] Denoting the tolerance by "tol", find $M, 0< M\leq L_1,$ such that
$$\sigma_1\geq \sigma_2\geq\cdots \geq\sigma_M>tol>\sigma_{M+1}\approx\cdots\approx\sigma_{L_1}\approx0.$$
If $M=L_1$, then increase $L_1$ and return to (2.a). \\
Otherwise, take the first $M$ columns of the matrix $V$ denoted by $V_0=V(:,1:M)$. Similarly, $W_0=W(:,1:M)$, and $\Sigma_0=diag(\sigma_1,\sigma_2,\cdots,\sigma_M)$.
    \item[(2.e)] Compute the eigenvalues $k_i$'s of $D=V_0^HC_{1,N_t}W_0\Sigma_0^{-1}\in\mathcal{C}^{M,M}$.
    \end{itemize}
\item[(3)] Validation. Compute the smallest eigenvalue $\lambda^0_i$ of $\mathbb{G}(k_i)$.
\begin{itemize}
    \item[(3.a)] Output $k_i$'s as eigenvalues if $|\lambda^0_i| < 10^{-12}$.
    \item[(3.b)] If $10^{-12} \le |\lambda^0_i| < 10^{-5}$, cover $R_i$ using smaller disks and go to Step (2).
    \item[(3.c)] If $|\lambda^0_i| \ge 10^{-5}$, discard $k_i$.
\end{itemize}
\end{itemize}


If one knows a priori a small region containing a few eigenvalues, one can skip Step (1) of {\bf Algorithm 1} and start Step (2) directly.
We summarize some guidelines when using the above algorithm in practice.

\begin{itemize}
    \item[{R1}:] Use a small tolerance in Step (2.d), e.g. tol$=1e-10$, for robustness.
    \item[{R2}:] Cover $R$ with smaller disks in Step (1) when possible.
    \item[{R3}:] Avoid the singularities of $\epsilon$, e.g., Drude-Sommerfel model.
\end{itemize}

\begin{remark}
The coefficients $\zeta_n(k)=\sqrt{k^2-\kappa_n^2}$ in the DtN map attain branch cuts, hence in the implementation of the algorithm, the region $R_i$ should not include the values $k=\kappa_n$, where the DtN map is not analytic.
\end{remark}

\section{Numerical examples}\label{NE}

\begin{figure}[!htbp]
 \begin{center}
 \hspace*{-2.5cm}
 \subfigure[Conducting metal]{\includegraphics[width=0.4\textwidth]{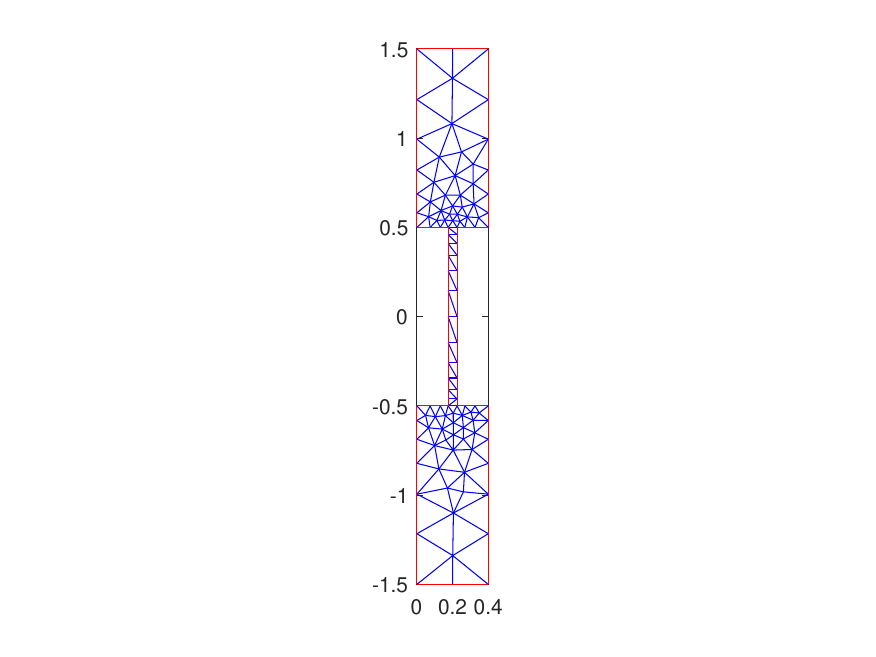}} \hspace*{-2cm}
 \subfigure[Real metal: rectangular slit]{\includegraphics[width=0.4\textwidth]{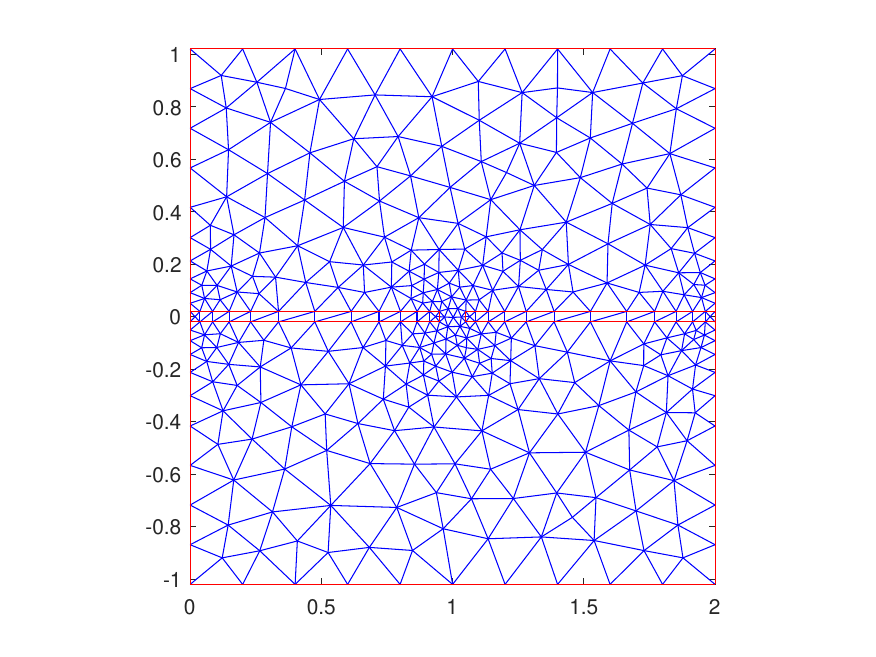}} \hspace*{-1cm}
 \subfigure[Real metal: trapezoidal slit]{\includegraphics[width=0.4\textwidth]{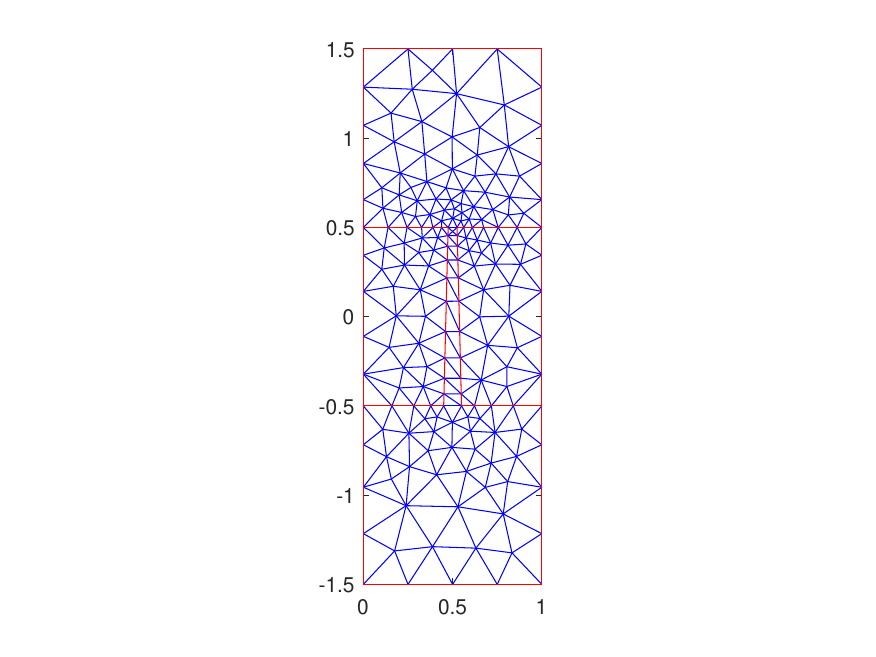}}
 \caption{Initial mesh $\mathcal T_{h_0}$ in the computational domain. The mesh is uniformly refined to obtain a series of mesh $\{\mathcal T_{h_j} \}_{j=0}^5$ for discretizing \eqref{eq:weak_main}. }\label{pec_example}
 \end{center}
 \end{figure}

In this section, we present several examples by considering different shapes of slit holes and electric permittivity functions. We first generate an unstructured initial mesh $\mathcal T_{h_0}$ for the computational domain ($\Omega_H\backslash\Omega$ for perfectly conducting metals or $\Omega_H$ for real metals), which is finer for small scale components of the domain and around the corners (see Figure \ref{pec_example}). Then the mesh is uniformly refined to obtain a series of meshes $\{\mathcal T_{h_j} \}_{j=0}^5$ and the linear Lagrange element is used for discretization to obtain \eqref{Gkv}. 
For the rest of this section, we call eigenvalues of \eqref{Gkv} small if their absolute values are small. Four examples are considered: (1) perfectly conducting metals; (2) Drude model without loss for the metal permittvity and rectangular slits; (3) Drude-Sommerfeld free electron model and rectangular slits; (4) Drude-Sommerfeld free electron model and trapezoidal slits.

\subsection{Perfectly conducting metal}


We first consider perfectly conducting metals  where the Neumann boundary condition is imposed over the metal boundary and the computational domain is $\Omega_H\backslash\Omega$.
For this configuration, the asymptotic expansions of eigenvalues for \eqref{Gkv} are available.
Assume that the period is $d$ and the thickness of the metallic slab is $\ell = 1$. The slit $S_0$ is a rectangle with width $\delta$.
The eigenvalues have the following asymptotic expansions for each $\kappa\in[-\frac{\pi}{d},\frac{\pi}{d}]$ (cf. \cite{linzha18a}):
\begin{equation}\label{kmkappa}
k_m(\kappa, \delta) = m\pi+2m\pi\left[\frac{1}{\pi}\delta \ln \delta+\left(\frac{1}{\alpha}+\gamma(m\pi,\kappa,d)\right)\delta\right]+O(\delta^2\ln^2\delta),\ \ \ m=1,2,3,\cdots,
\end{equation}
where the constant $\alpha\approx-1.1070218960566$ and
\begin{equation}
\gamma(k,\kappa,d)=\frac{1}{\pi}\left(3\ln 2+ln\frac{\pi}{d}\right)-\frac{i}{d}\frac{1}{\zeta_0(k)}+\sum_{n\neq0}\left(\frac{1}{2\pi|n|}-\frac{i}{d}\frac{1}{\zeta_n(k)}\right).
\end{equation}

Let $d = 0.4$ and the Bloch wavenumber $\kappa=\frac{\pi}{d}$.
By neglecting the high-order term $O(\delta^2ln^2\delta)$ in the asymptotic expansion \eqref{kmkappa}, the smallest eigenvalues for $\delta=0.05,0.02,0.01$ are given by
\[
k_1 \left(\frac{\pi}{d}, 0.05\right) \approx 2.8146, \quad k_1 \left(\frac{\pi}{d}, 0.02\right) \approx 2.9741, \quad k_1 \left(\frac{\pi}{d}, 0.01\right) \approx 3.0440.
\]

\begin{figure}[!htbp]
\centering
\includegraphics[width=0.7\textwidth]{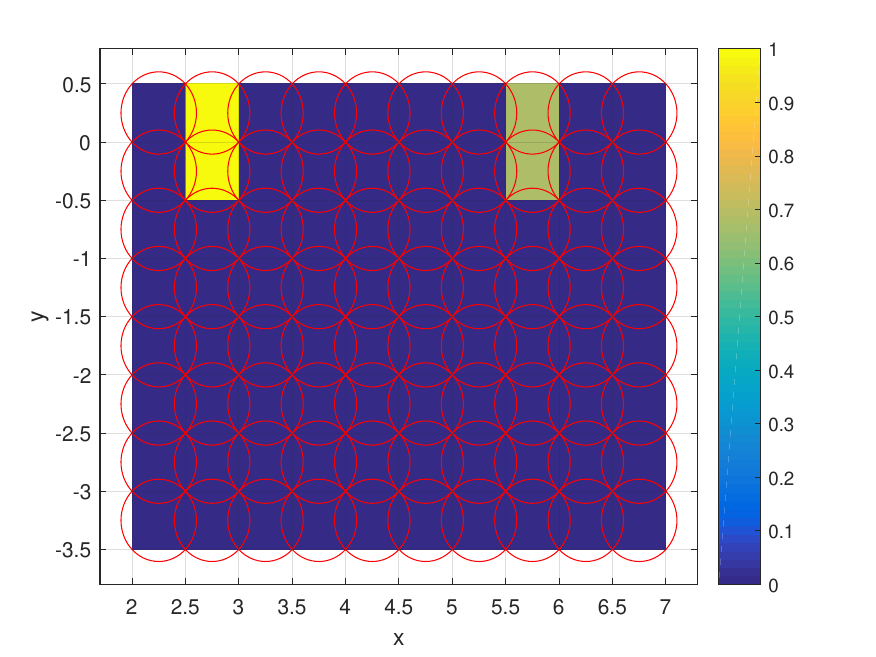}
\caption{Indicators of sub-regions in Step (1) of Algorithm 1.}
\label{fig:indicator_step1}
\end{figure}

We first demonstrate the process of Step (1) in Algorithm 1. Set the slit width $\delta=0.05$ and use a mesh with $h\approx0.025$. Assume the search region is $R=[2,7]\times[-3.5,0.5]$ in the fourth quadrant of $\mathbb C$. We use $80$ uniform disks to cover $R$.
The normalized indicators $\|P_i v\|$ are shown in Figure \ref{fig:indicator_step1}. There are four disks with large indicators that are kept for Step (2).

Next we check the convergence of the smallest eigenvalues with respect to the mesh size. The initial mesh size is denoted by $h_0$ and $h_j = \frac{h_{j-1}}{2}$ $(j = 1, 2,\cdots, 5)$ for the subsequent refined meshes.
The relative convergence order is defined as
$$
\log_2\left(\left|\frac{k^j-k^{j-1}}{k^{j+1}-k^j}\right|\right),\ \ \ j=1,2,3,4,
$$
where $k^j$ is the smallest eigenvalue computed using ${\mathcal T}_{h_j}$.

The DOFs (degrees of freedoms) for $\delta=0.05,0.02,0.01$ on the finest meshes are $83329$, $161665$, $229505$,  respectively.
For the DtN mapping, we take $D_t=50$. The number of equidistant nodes on $\partial R$ is $N_t=64$.
The computed eigenvalues are shown in Table \ref{table4}. The eigenvalues converge as the mesh is refined and the convergence rate is less than $2$.

\begin{table}
\begin{center}
\caption{The computed eigenvalues for the perfectly conducting metal with rectangular slits of width $\delta$. The search domain $R$ is the disk centered at $(3,0)$ of radius $r=0.5$.}
\medskip
\label{table4}
\begin{tabular}{ccccccc}
\hline
Mesh refinement &$\delta=0.05$&Order& $\delta=0.02$&Order&$\delta=0.01$&Order\\
\hline
0   &2.91449717&       &3.01693273&       &3.07038407&    \\
1   &2.87515496&       &2.99555023&       &3.05464236&    \\
2   &2.86020713&1.3961 &2.98803899&1.5093 &3.04969882&1.6701\\
3   &2.85449203&1.3871 &2.98533567&1.4743 &3.04809664&1.6255\\
4   &2.85229502&1.3792 &2.98433941&1.4401 &3.04755690&1.5679\\
5   &2.85144330&1.3671 &2.98396374&1.4071 &3.04736774&1.5127\\
\hline
\end{tabular}
\end{center}
\end{table}


\subsection{Sheetmetal grating with rectangular slits}
We consider the sheetmetal grating with rectangular slits and compare the result with that in \cite{Lin}, where a mode matching method is applied. The parameters used for the metallic grating are: $\ell=40$nm, $d=2\mu$m, and $\delta=0.1\mu$m. The permittivity of the metal is given by the Drude model without loss:
\[
\varepsilon_m(\omega) =  \varepsilon_0 \left( 1- \frac{\omega_p^2}{\omega^2}\right),
\]
where $\varepsilon_0$ is the permittivity in the vacuum and $\omega_p$=300THZ is the plasma frequency.

We employ a scaling for the geometry with a factor of $\alpha=10^6$ ($1\mu$m to $1$m) such that $\ell = 0.04$m, $d=2$m and $\delta=0.1$m.
The wavenumber becomes $\hat k = k/\alpha$ in \eqref{eq:weak_main} and $\varepsilon_m$ can be written as
$$ \varepsilon_m(\omega) =  \varepsilon_0 \left( 1- \frac{\omega_p^2}{\omega^2}\right) = \varepsilon_0 \left( 1- \frac{\omega_p^2}{(ck)^2}\right) =
\varepsilon_0 \left( 1- \frac{\omega_p^2}{(c\alpha \hat k)^2}\right) = \varepsilon_0 \left( 1- \frac{\hat\omega_p^2}{\hat k^2}\right), $$
where $c=3\times 10^8$m/s is the speed of light and $\hat\omega_p = \omega_p/(c\alpha)$ is the scaled plasma frequency.
Note that the frequency $\omega=ck=c\alpha\hat{k}$ in \cite{Lin}.

\begin{figure}[!htbp]
\begin{center}
\subfigure[]{\includegraphics[width=0.45\textwidth]{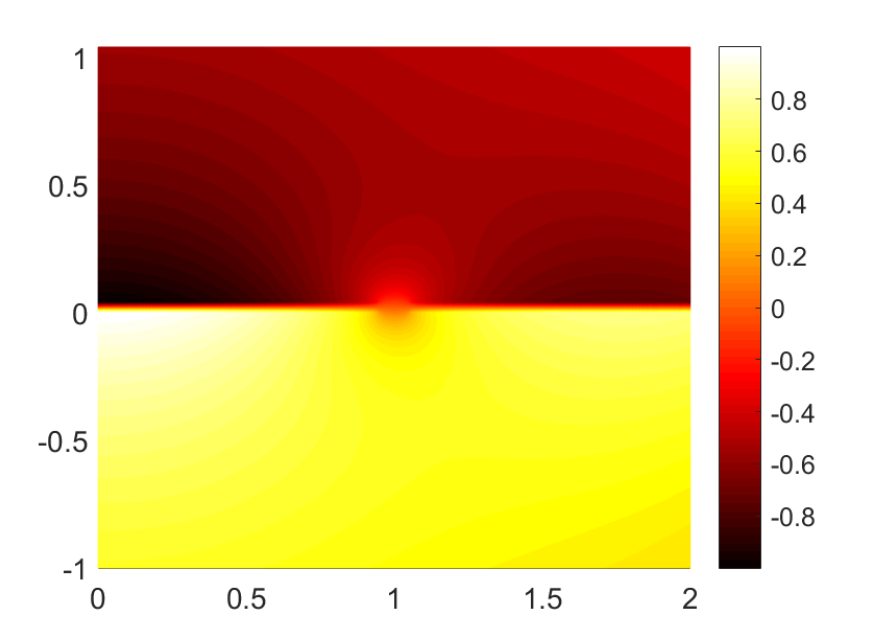}}
\subfigure[]{\includegraphics[width=0.45\textwidth]{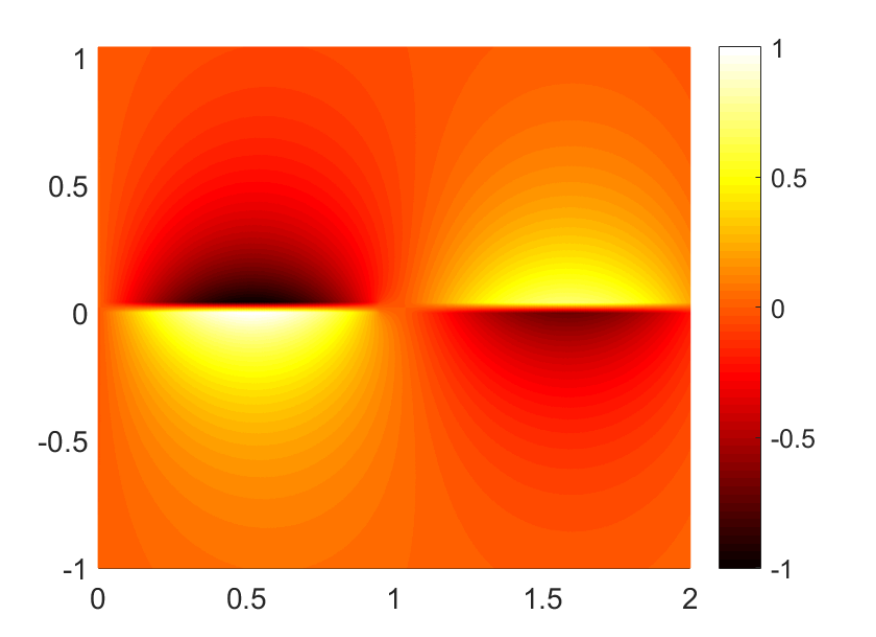}}\\
\subfigure[]{\includegraphics[width=0.45\textwidth]{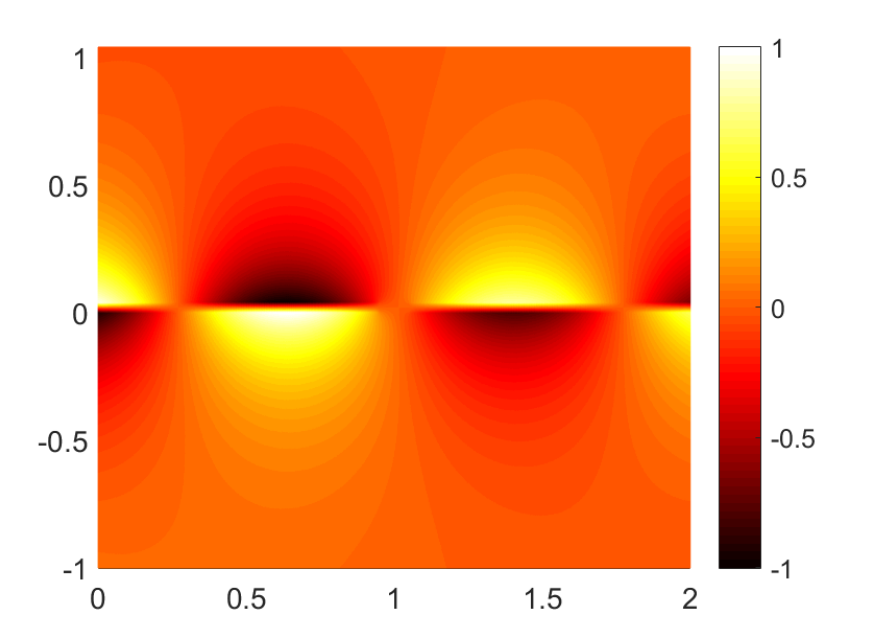}}
\subfigure[]{\includegraphics[width=0.45\textwidth]{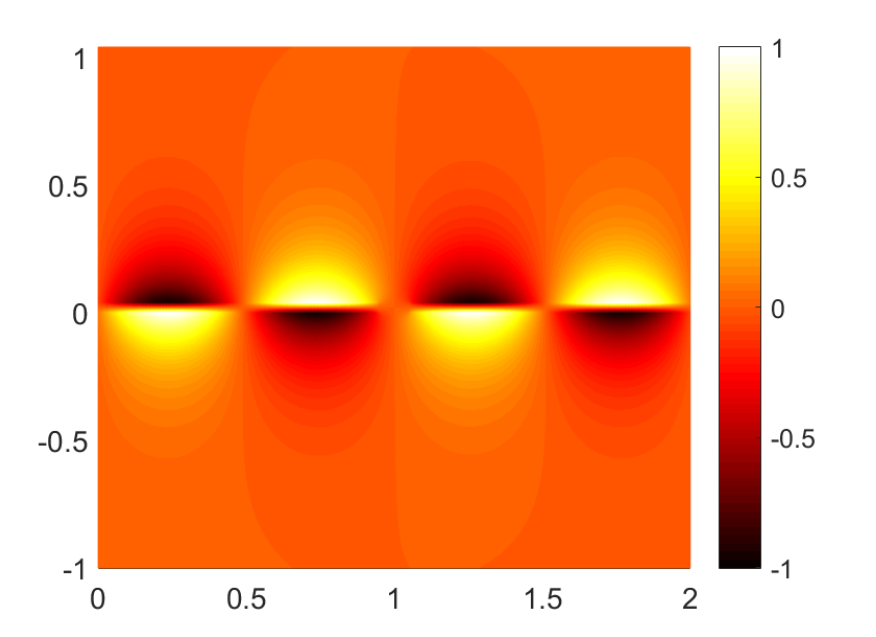}}
\caption{Real parts of eigenfunctions when the Bloch wavenumber
$\kappa=\frac{\pi}{4d}$ for the sheetmetal grating with rectangular slits.}\label{newexample_kappa1}
\end{center}
\end{figure}

The initial mesh is shown in Figure \ref{pec_example} (b). We refine the initial mesh 4 times and end up with 87872 DOFs. The truncation order for the DtN map is $D_t=50$. The initial search region is a disk centered at $(0.2,0)$ with radius $r=0.15$. For $\kappa = \frac{\pi}{4d}$, $4$ eigenvalues are obtained:
\begin{center}
$\hat{k}_1$ = 0.12492920,\ \ \ $\hat{k}_2$ = 0.23916592,\ \ \
$\hat{k}_3$ = 0.27838236,\ \ \ $\hat{k}_4$ = 0.33281163.
\end{center}
The corresponding frequencies are (in THz)
$$\omega_1=37.478757,\quad \omega_2=71.749776,\quad \omega_3=83.514708,\quad \omega_4=99.843492.$$
The associated eigenfunctions are shown in Figure~\ref{newexample_kappa1}.

\begin{figure}
\centering
\includegraphics[width=0.5\textwidth]{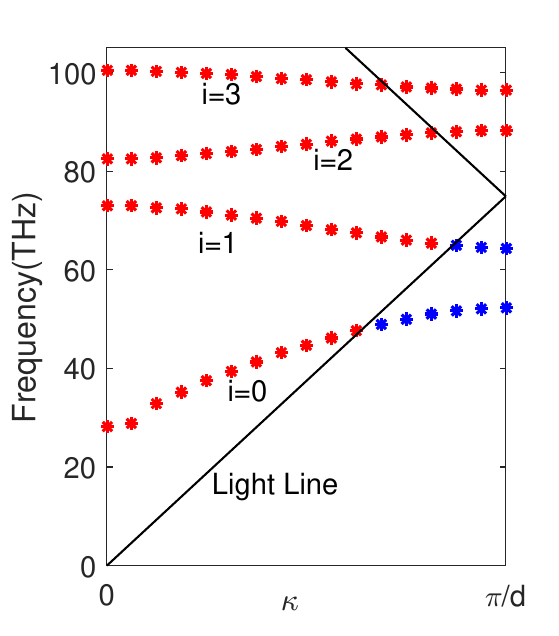}
\caption{The band structure for the sheetmetal grating, where the Drude model \eqref{FreeElectron} is used for the permittivity of the metal.}
\label{fig:band_structure666666}
\end{figure}

Figure \ref{fig:band_structure666666} shows the band structure for $\kappa\in[0,\pi/d]$. The dispersion curves indexed by $i=1,2,3$ are consistent with \cite{Lin}. However, the dispersion of unbounded SPP modes (red dots) for $i=0$ is different than that in Figure 2 of \cite{Lin}. We obtain a continuous dispersion curve throughout the Brillouin zone, while the dispersion curve in \cite{Lin} is not continuous. Such discrepancy is subjected to further investigation.

%

\subsection{Metallic grating with Drude-Sommerfeld model}
We consider a metallic grating with thickness $\ell$ and width $\delta$ for the rectangular slit holes.  The Drude-Sommerfeld model (\ref{eq:em}) is used for the metal permittivity
\begin{equation}\label{FreeElectron}
\varepsilon_{m}(\omega)=1-\frac{\omega_p^2}{\omega^2+\Gamma ^2}+i\frac{\Gamma\omega_p^2}{\omega(\omega^2+\Gamma^2)}.
\end{equation}
The plasma frequency and the damping constant for gold are $\omega_p=1.38\times 10^{16}/s$ and $\Gamma = 1.075\times 10^{14}/s$.
We consider a multiscale structure by assuming that the metal thickness $\ell=0.1\mu$m, the period $d=0.1\mu$m, and the width $\delta=5$nm for the rectangular slit.

Using a scaling factor $\alpha_0=10^7$ and denoting the frequency by $\omega=ck=c\alpha_0\hat{k}$, the permittivity can be written as
\begin{equation}
\varepsilon_{m}(\omega)=1-\frac{\hat{\omega}_p^2}{\hat{k}^2+\hat{\Gamma} ^2}+i\frac{\hat{\Gamma}\hat{\omega}_p^2}{\hat{k}(\hat{k}^2+\hat{\Gamma}^2)},
\end{equation}
where $\hat{\omega}_p=\frac{\omega_p}{c\alpha_0}$ and $\hat{\Gamma}=\frac{\Gamma}{c\alpha_0}$.
The scaled geometry parameters are $\hat{\ell}=1m$, $\hat{d}=1m$ and $\hat{\delta}=0.05m$. We set the Bloch wavenumber $\kappa = \frac{\pi}{2}$ and $D_t=100$. The mesh is much finer around the slit (see Figure~\ref{pec_example}).

We consider several search regions for this example.
The first search region $R_1$ is the disk centered at $(2.4,0)$ with radius $0.3$.
Step (1) of \textbf{Algorithm 1} indicates that $R_1$ contains eigenvalues.
Step (2) output $7$ values $\Lambda_i, i=1, \ldots 7$. We list them and the smallest eigenvalues $\lambda_i^0$'s of $\mathbb{G}(\Lambda_i)$ below
\[
\begin{array}{ll}
  2.24459845 - 0.01762724i,& 10^{-14}\cdot(-2.30552126 + 2.43361749i),\\
  2.39786005 - 0.01791077i,& 10^{-15}\cdot( 7.75373810 - 8.05338438i),\\
  2.51904653 - 0.01788382i,& 10^{-14}\cdot(-2.78954512 - 5.84148363i),\\
  2.61689989 - 0.01772865i,& 10^{-15}\cdot(-1.64818948 - 2.31104652i),\\
  2.66286367 - 0.01832286i,& 10^{-15}\cdot( 1.79896380 + 4.00565464i),\\
  2.66797279 - 0.01831849i,& 10^{-16}\cdot( 1.06438139 + 3.72171976i),\\
  2.69949977 - 0.01779482i,& 10^{-15}\cdot(-9.38044522 + 3.79717998i).
\end{array}
\]
Since the smallest eigenvalues $\lambda_i^0$'s are smaller than $10^{-12}$, all $\Lambda_i, i=1, \ldots, 7,$ are eigenvalues of $\mathbb{G}(\cdot)$.

The second region $R_2$ is the disk centered at $(1.5,0)$ with radius $r=0.7$. Again, Step (1) of \textbf{Algorithm 1} indicates that $R_2$ is a region containing eigenvalues. Step (2) outputs $12$ values $\Lambda_1,\Lambda_2,\cdots, \Lambda_{12}$,
\[
\begin{array}{rr}
\Lambda_1=0.81562725 - 0.01021274i, &\Lambda_2=1.38497728 - 0.01161541i,\\
\Lambda_3=1.50732390 - 0.00098266i, &\Lambda_4=1.51821059 - 0.00433577i,\\
\Lambda_5=1.55863415 - 0.16599511i, &\Lambda_6=1.56182180 - 0.16708495i,\\
\Lambda_7=1.61670231 - 0.42231762i, &\Lambda_8=1.62251917 - 0.42252083i,\\
\Lambda_9=1.68882271 - 0.61675923i, &\Lambda_{10}=1.69054224 - 0.61664504i,\\
\Lambda_{11}=1.78400381 - 0.01661049i, &\Lambda_{12}=2.04682050 - 0.01695064i.
\end{array}
\]

Plugging these values into $\mathbb{G}(\cdot)$ for Step (3) , the smallest eigenvalues $\lambda_{\Lambda_i}^0$'s of $\mathbb{G}(\Lambda_i)$ are as follows
\[
\begin{array}{rl}
     10^{-11} \cdot (4.87892152+3.44068378i),& 10^{-8}\cdot (-5.80955512 - 2.92448205i), \\
     10^{-6} \cdot (1.64592598 - 0.39010526i), & 10^{-6}\cdot(-1.41048364 + 0.61461716i), \\
    10^{-5}\cdot (0.04797358 + 4.05544658i), & 10^{-5}\cdot(-0.00236747 + 4.09881869i).\\
    10^{-5}\cdot (0.54966602 + 8.49807850i),& 10^{-5}\cdot(0.48749558 + 8.52658665i),\\
 10^{-5}\cdot (-6.14003880 + 2.21285125i),& 10^{-5}\cdot(-6.15257226 + 2.21529267i),\\
 10^{-10}\cdot(0.10128613) - 1.31158998i),& 10^{-12} \cdot (-3.20743760+2.61646040i).
\end{array}
\]
Although none of them is smaller than $10^{-12}$ in norm, but some values are small enough for further investigation. We use disks centered at these values with smaller radius, e.g., $r=0.03$, as the input regions for Step (2) of {\bf Algorithm 1}. We find $6$ of them are eigenvalues, which are listed below together with the corresponding smallest eigenvalues of $\mathbb{G}(\Lambda_i)$:
\[
\begin{array}{ll}
0.81562760 - 0.01021249i, & 10^{-17}\cdot(3.20096386 + 0.06686478i), \\
1.38472724 - 0.01173458i, & 10^{-17} \cdot(-0.84617653 - 2.41196884i),\\
1.51324908 - 0.00230914i,& 10^{-17} \cdot(-7.40234269 + 1.53633700i),\\
1.52492183 - 0.00437709i,& 10^{-17} \cdot(1.08998605  + 5.58769924i),\\
1.78400511 - 0.01661062i ,& 10^{-17} \cdot(-8.21478417 - 0.00749054i),\\
2.04682049 - 0.01695064i ,& 10^{-14} \cdot(-1.86190705 - 1.67313875i).
\end{array}
\]


Other values that are discarded, e.g., $\Lambda_5$. We show how this is decided. When the search region $R_3$ is the disc centered at  $\Lambda_5$ with radius $0.03$, Step 2 computes $6$ values inside $R_3$. We list them and the corresponding smallest eigenvalues of $\mathbb{G}(\cdot)$ below
\[
\begin{array}{ll}
1.57869531 - 0.14731601i,& 10^{-5} \cdot (-4.69518900 - 0.63047006i), \\
1.57892533 - 0.14727291i,& 10^{-5} \cdot (-4.69583811 - 0.63208982i),\\
1.58043466 - 0.16339190i,& 10^{-5} \cdot (-4.79919586 - 0.57366528i),\\
1.58100471 - 0.16323098i,& 10^{-5} \cdot (-4.80075913 - 0.57765327i),\\
1.58188508 - 0.18004639i,& 10^{-5} \cdot (-4.89733368 - 0.50639969i),\\
1.58210712 - 0.17996384i,& 10^{-5} \cdot (-4.89793654 - 0.50795186i).
\end{array}
\]
Since the smallest eigenvalue of $\mathbb{G}(\Lambda)$ is of magnitude $10^{-5}$, $\Lambda_5$ is discarded. In fact, if one continues to check a series of smaller regions, no eigenvalues will be found.

\begin{figure}
\centering
\includegraphics[width=0.6\textwidth]{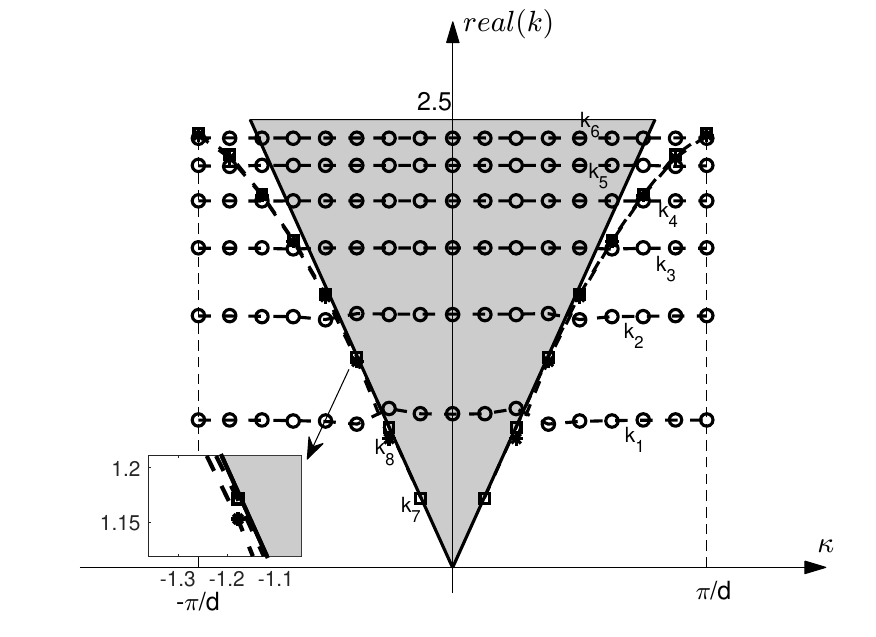}
\caption{The band structure for the metallic grating considered in Section 5.3. The X-axis and Y-axis represents the Bloch wavenumber $\kappa$ and the real part of resonances respectively.}
\label{fig:band_structure22}
\end{figure}
\begin{figure}
\centering
\subfigure[\textbf{$u(k_1;\cdot)$}]{\includegraphics[width=0.3\textwidth]{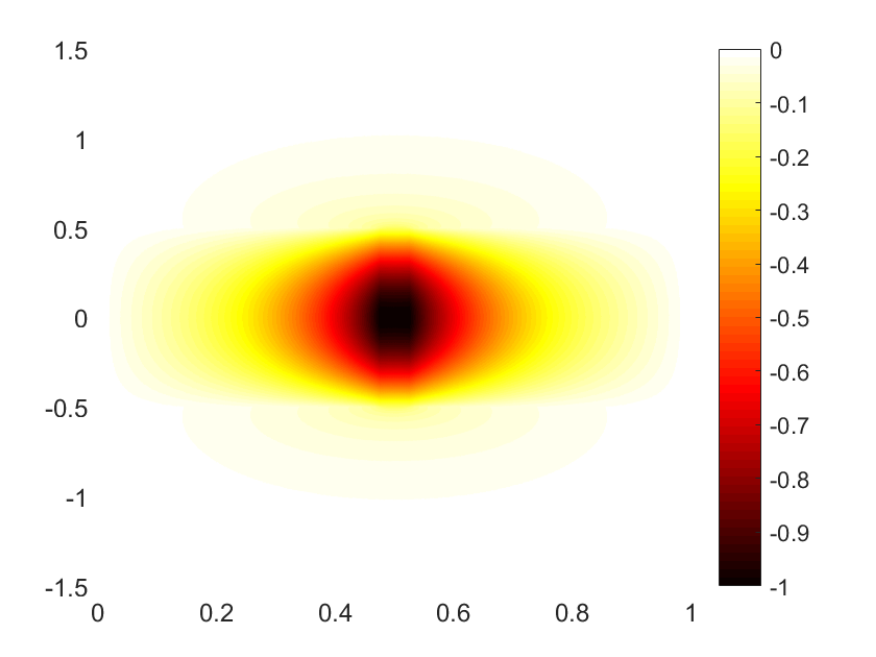}}
\subfigure[\textbf{$u(k_2;\cdot)$}]{\includegraphics[width=0.3\textwidth]{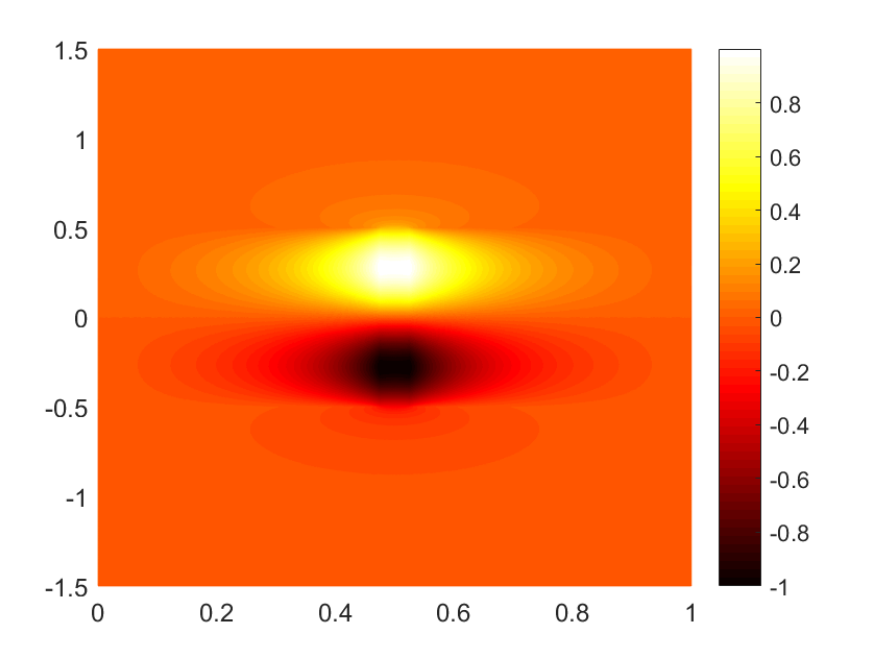}}
\subfigure[\textbf{$u(k_3;\cdot)$}]{\includegraphics[width=0.3\textwidth]{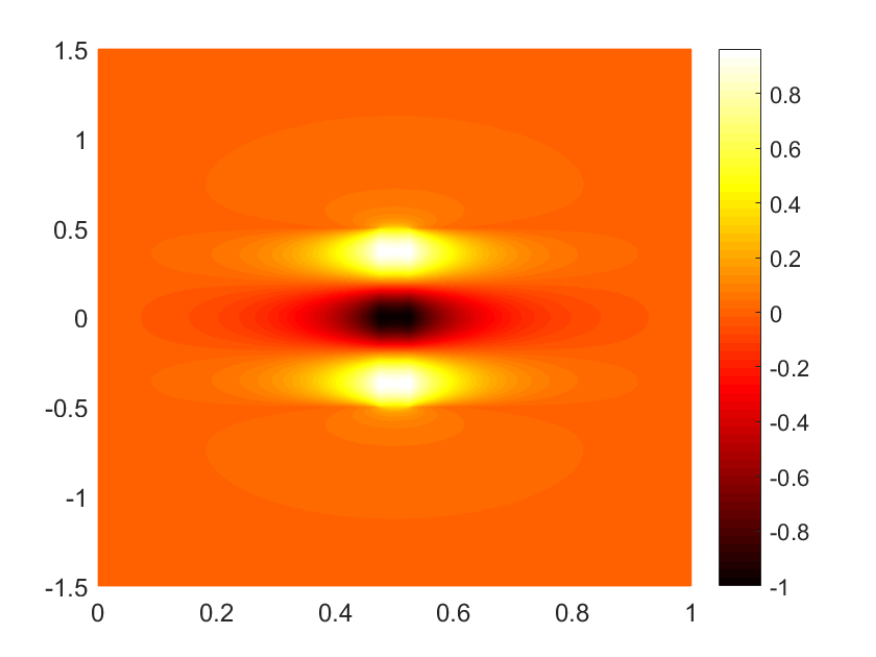}}
\subfigure[\textbf{$u(k_4;\cdot)$}]{\includegraphics[width=0.3\textwidth]{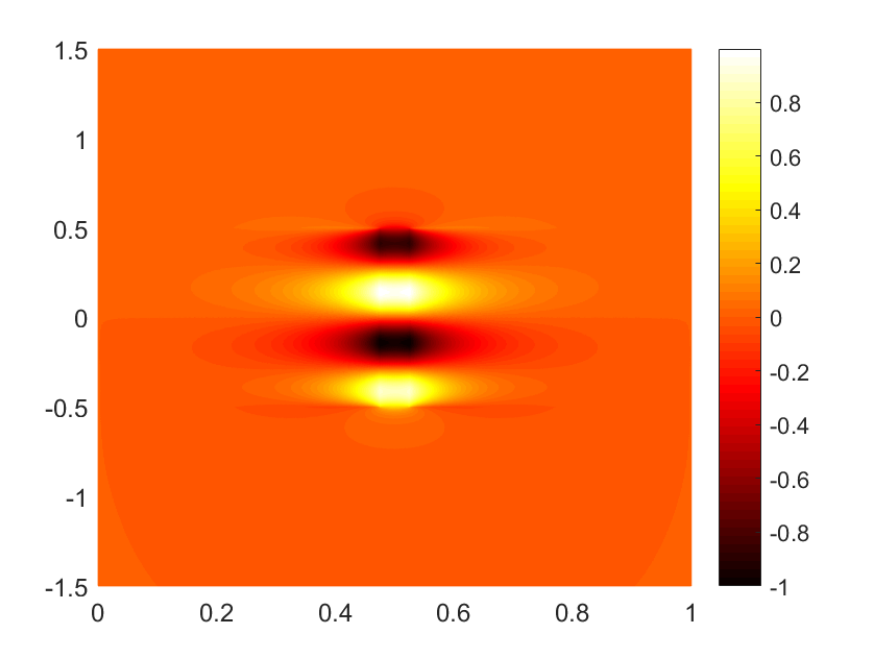}}
\subfigure[\textbf{$u(k_5;\cdot)$}]{\includegraphics[width=0.3\textwidth]{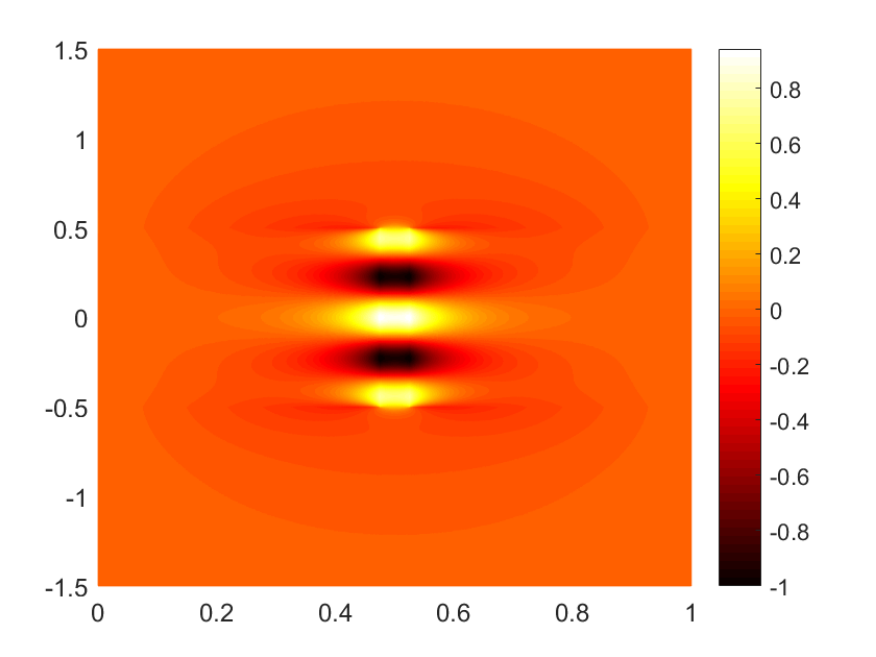}}
\subfigure[\textbf{$u(k_6;\cdot)$}]{\includegraphics[width=0.3\textwidth]{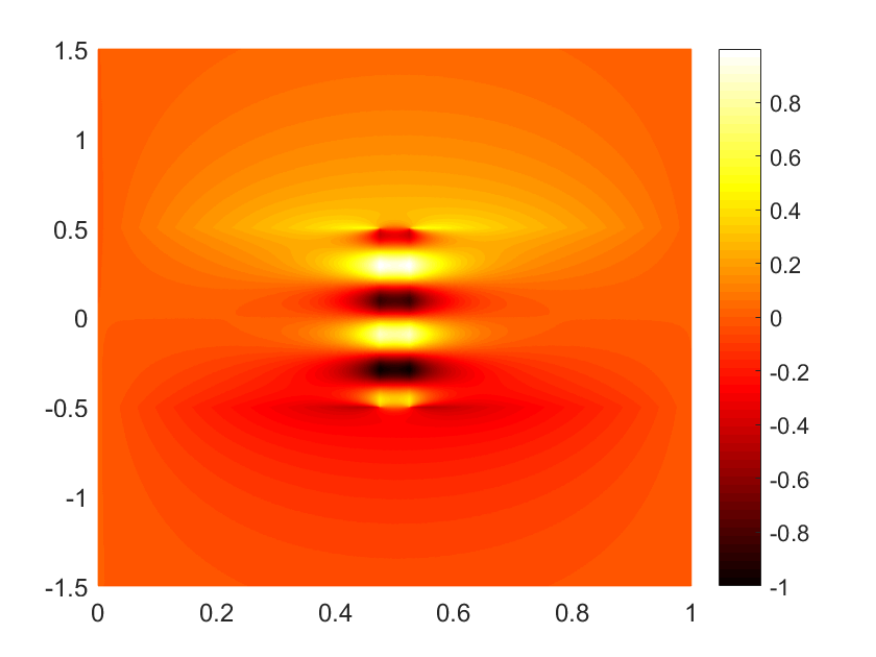}}
\subfigure[\textbf{$u(k_7;\cdot)$}]{\includegraphics[width=0.3\textwidth]{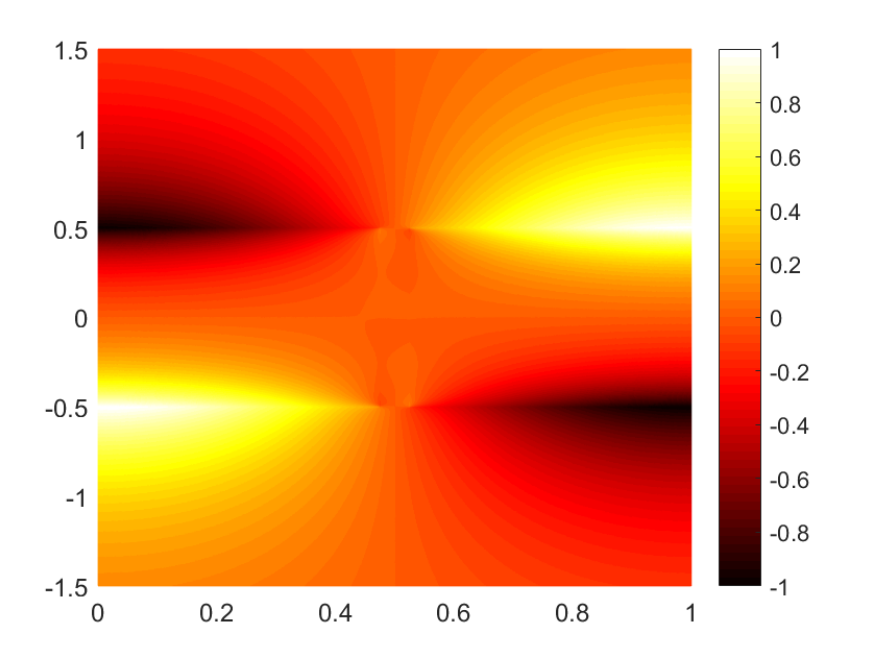}}
\subfigure[\textbf{$u(k_8;\cdot)$}]{\includegraphics[width=0.3\textwidth]{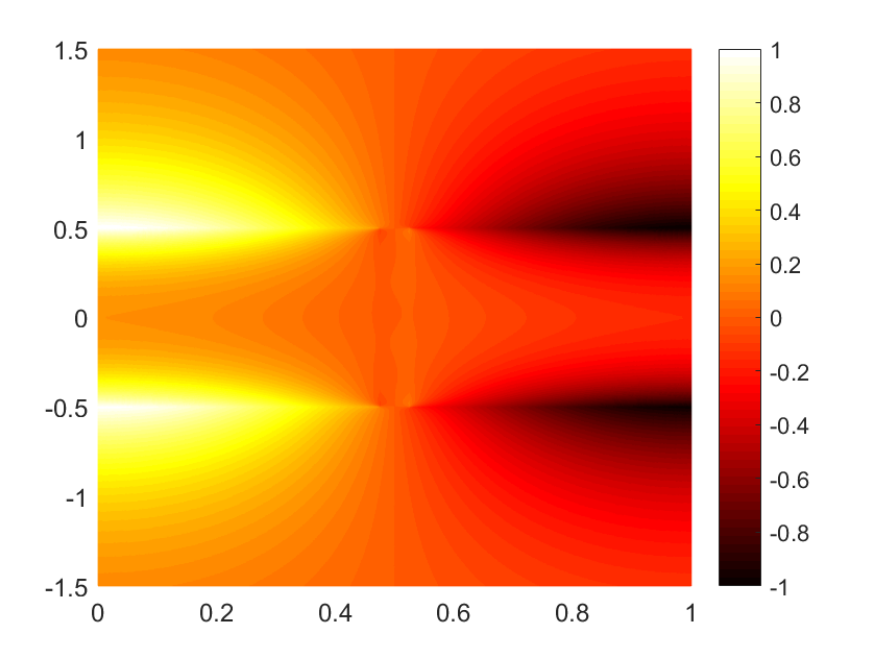}}
\caption{Real parts of eigenfunctions at $\kappa=\pi$ for the metallic grating considered in Section 5.3.}
\label{eigenfunctions8888}
\end{figure}

The band structure of the metallic grating is shown in Figure \ref{fig:band_structure22}. The bands $k_j(\kappa)$ ($1\le j \le 6)$ are resonances induced by the slit holes and the bands $k_j(\kappa)$ ($7\le j \le 8)$ are surface plasmonic resonances. More specifically,
when $\kappa=\pi$, the computed resonances are
\begin{align}\nonumber
k_1&=0.82333707 - 0.01098713i,\quad k_2=1.40413513 - 0.01417461i,\\ \nonumber
k_3&=1.78249483 - 0.01600898i,\quad k_4=2.04659065 - 0.01685595i, \\ \nonumber
k_5& =2.24213036 - 0.01717807i,\quad k_6=2.38932484 - 0.01689001i,\\ \nonumber
k_7&=2.41320003 - 0.00952127i, \quad k_8=2.42594474 - 0.00918862i.
\end{align}

The corresponding eigenfunctions are shown in Figure~\ref{eigenfunctions8888}.
For resonances induced by the slit holes, the eigenfunctions are localized inside and around slit holes (see Figure \ref{eigenfunctions8888} (a)-(f)). In addition, they become more oscillatory for higher-order modes. We call these resonances cavity resonances, as they are induced by the cavities formed by the slit holes. Their studies for the perfect conducting metals are reported in \cite{linzha18a}. We see that for real metals considered here, the wave inside the slit holes penetrates into the metal with skin depth depending the metal loss, which shifts the resonance frequencies significantly.
On the other hand, for surface plasomic resonances, the eigenfunctions are localized at the metal surface (see Figure \ref{eigenfunctions8888} (g)-(h)). These eigenmodes harvest energy near the metal and are typical surface plasmonic modes induced by the metal \cite{maier07}.

\newpage

\subsection{Metallic grating with trapezoidal slit holes}
Finally, we examine a model when the slit hole is a trapezoid. The initial mesh for the discretization is shown in Figure \ref{pec_example}(c).
The finite element discretization allows more general geometry, which makes it attractive in practical computation. The permittivity for gold metal  is again the Drude-Sommerfeld model \eqref{FreeElectron}. Using a scaling factor $\alpha_0= 10^7$,
the metal thickness is $\ell=1$ and the period is $d=1$. The isosceles trapezoidal slit $S_0$ has a top width $0.05$ and a base width $0.1$ as shown in Figure \ref{pec_example} (c).

\begin{figure}
\centering
\subfigure[\textbf{$u(k_1;\cdot)$}]{\includegraphics[width=0.3\textwidth]{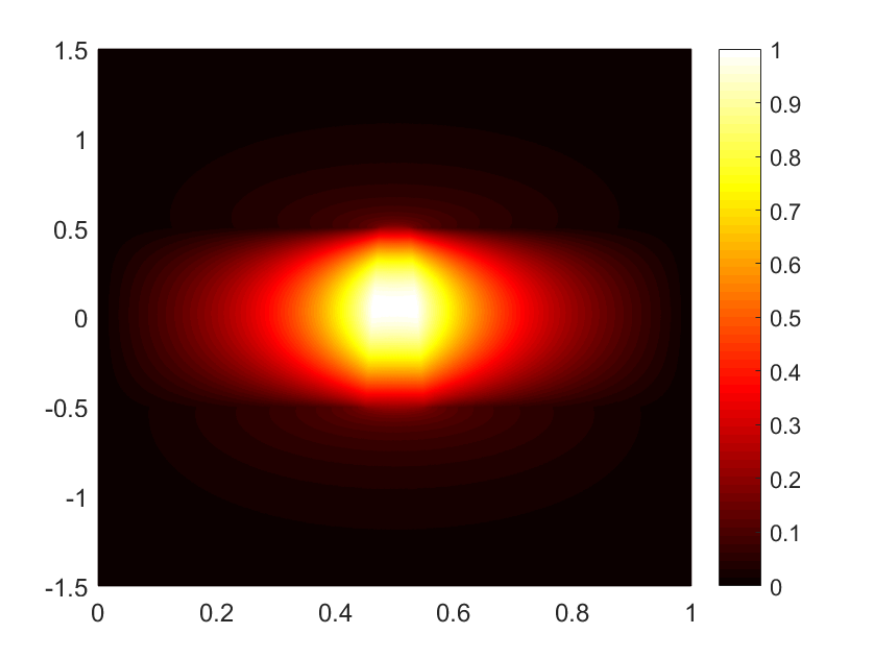}}
\subfigure[\textbf{$u(k_2;\cdot)$}]{\includegraphics[width=0.3\textwidth]{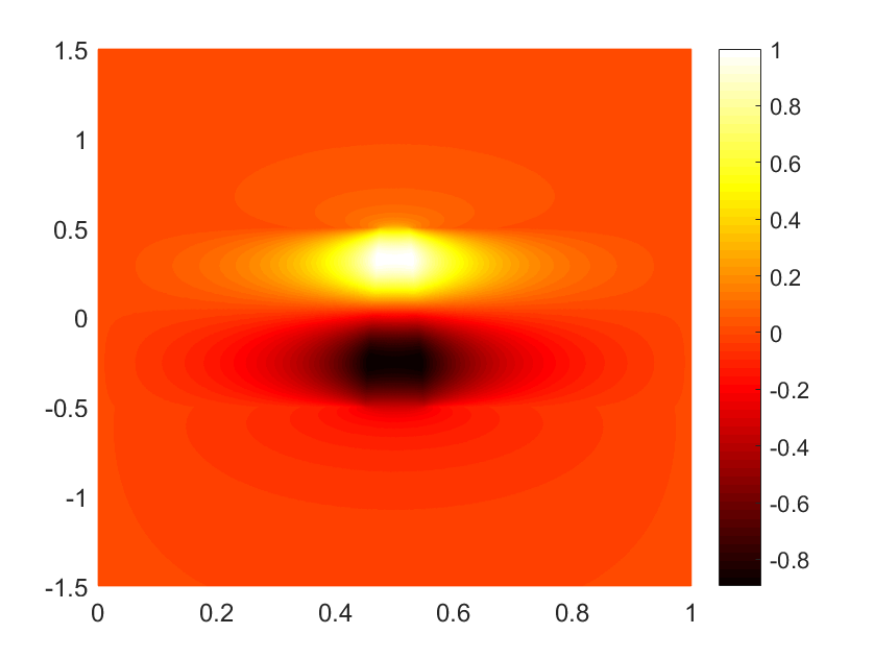}}
\subfigure[\textbf{$u(k_3;\cdot)$}]{\includegraphics[width=0.3\textwidth]{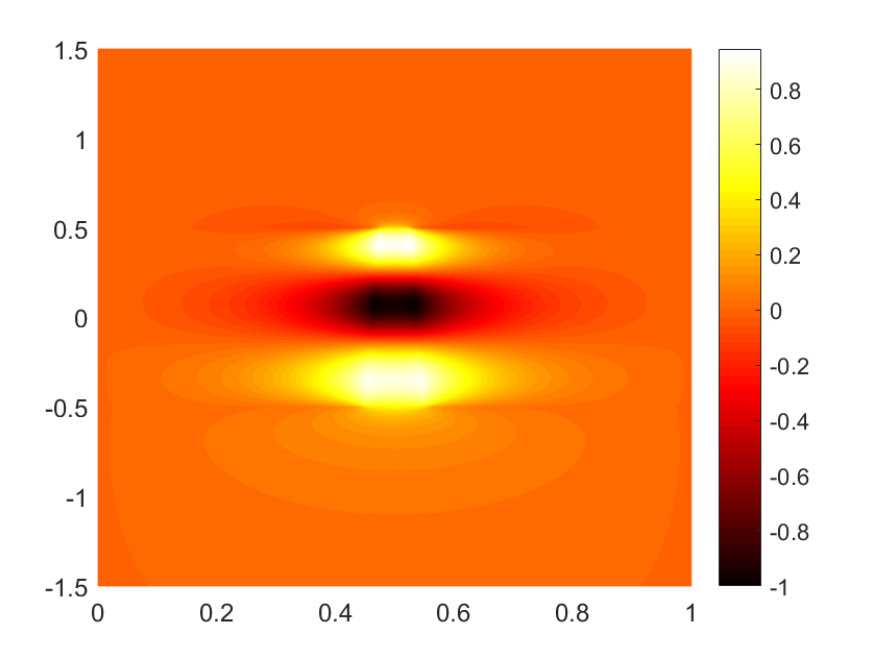}}
\subfigure[\textbf{$u(k_4;\cdot)$}]{\includegraphics[width=0.3\textwidth]{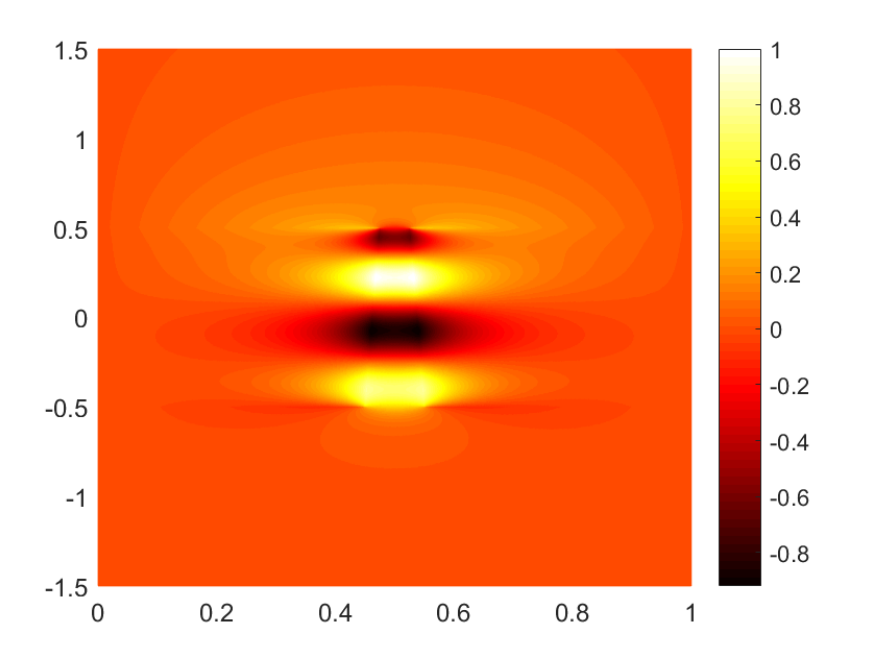}}
\subfigure[\textbf{$u(k_5;\cdot)$}]{\includegraphics[width=0.3\textwidth]{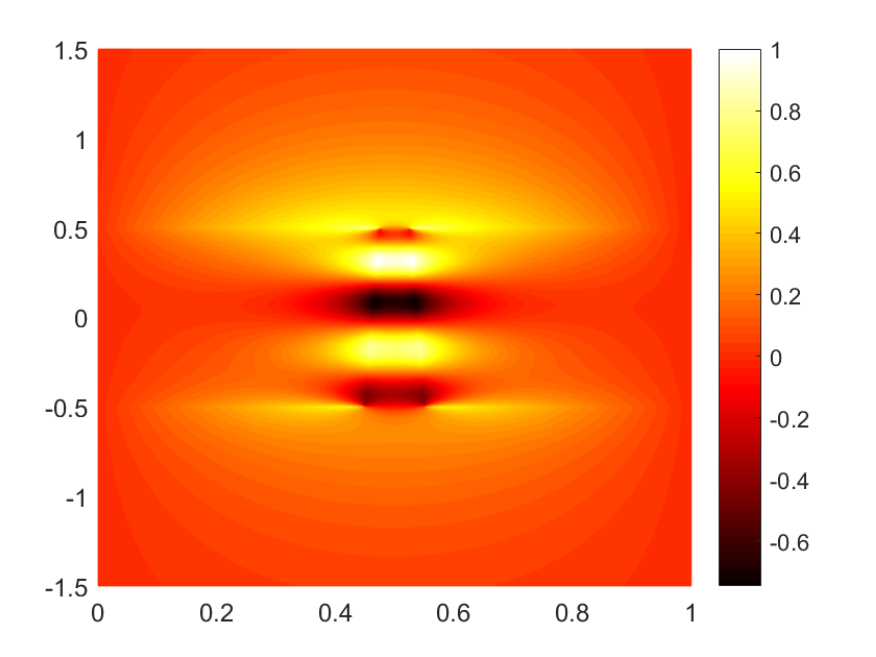}}
\subfigure[\textbf{$u(k_6;\cdot)$}]{\includegraphics[width=0.3\textwidth]{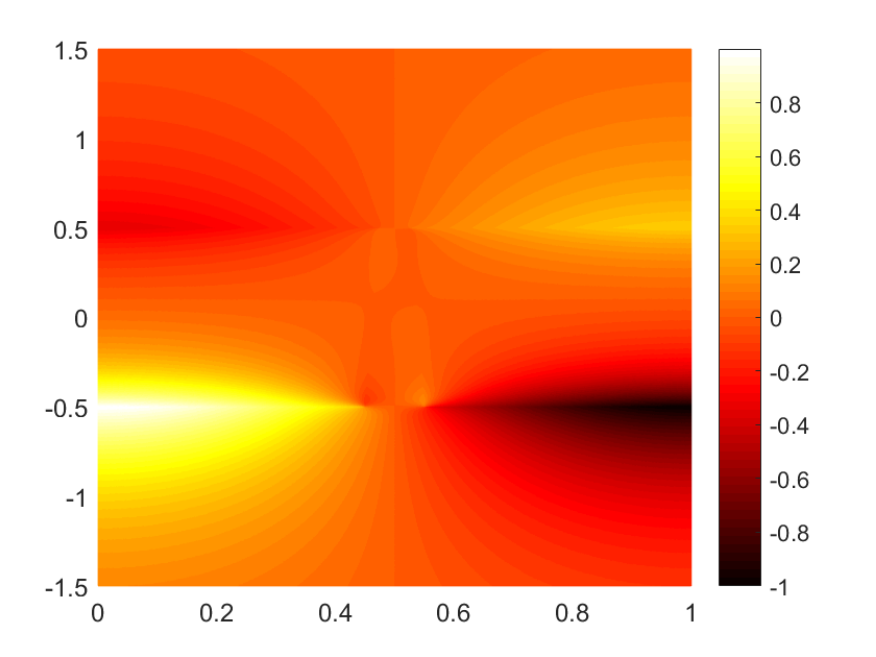}}
\subfigure[\textbf{$u(k_7;\cdot)$}]{\includegraphics[width=0.3\textwidth]{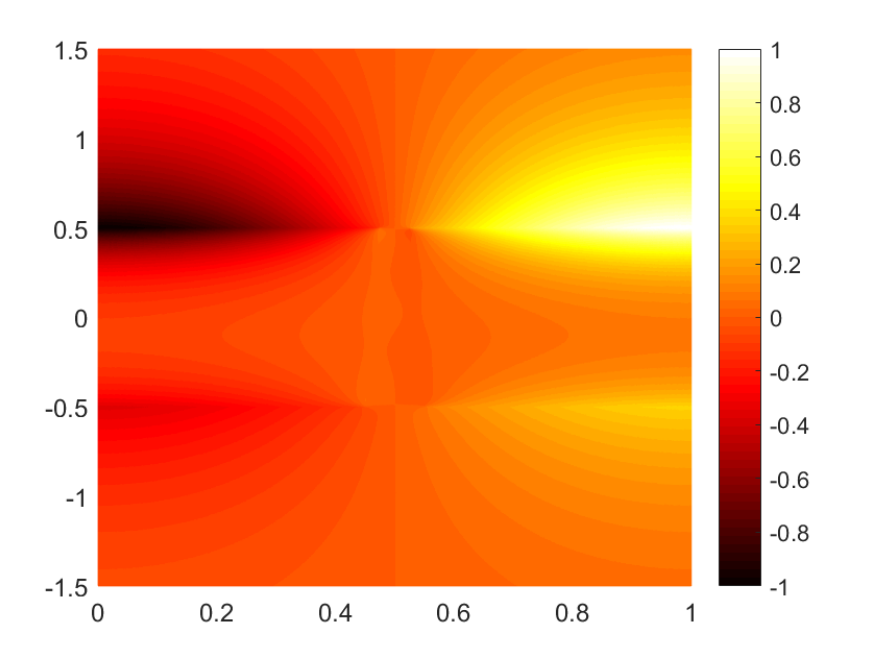}}
\caption{Real parts of eigenfunctions at $\kappa=\pi$ for the metallic grating considered in Section 5.4.}
\label{eigenfunctions9999}
\end{figure}

\begin{figure}
\centering
\includegraphics[width=0.6\textwidth]{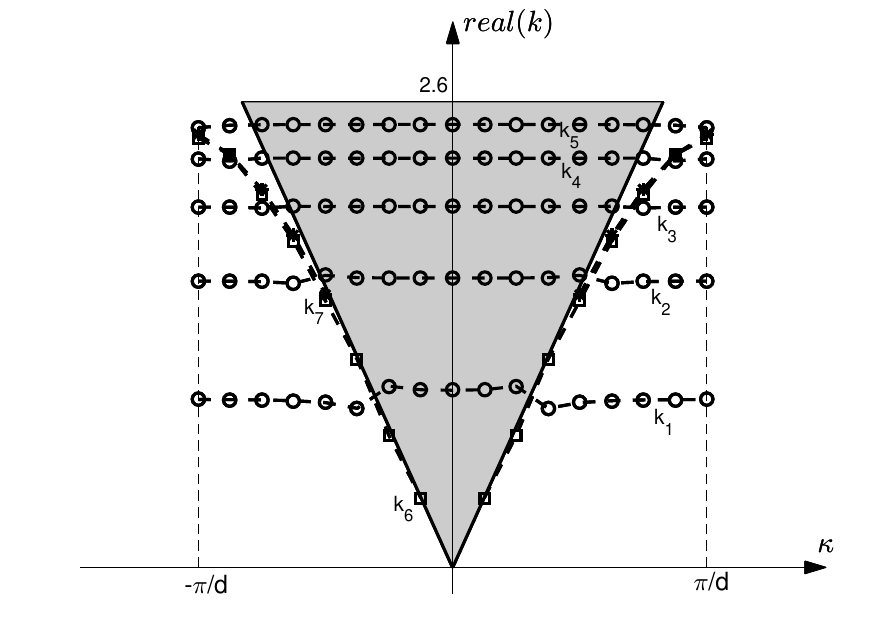}
\caption{The band structure for the metallic grating considered in Section 5.4.}
\label{fig:band_structure2}
\end{figure}

%
We use a mesh with a total of $43776$ DOFs and set $D_t=100$, $N_t=256$.
The search region $R_1$ is the disc centered at $(1.5,0)$ with radius $r=1$. 
The computed eigenvalues at $\kappa=\pi$ are
\begin{align}\nonumber
k_1&=0.93798322 - 0.01061232i, \quad k_2=1.59880173 - 0.01400498i,\\ \nonumber
k_3&=2.01182415 - 0.01597829i, \quad k_4=2.28107400 - 0.01677095i,\\ \nonumber
k_5&=2.45552736 - 0.01595151i, \quad k_6= 2.39889904 - 0.01005366i,\\ \nonumber
k_7&=2.41709341 - 0.00950237i.
\end{align}
The corresponding eigenfunctions are shown in Figure \ref{eigenfunctions9999}. The band structure for all eigenvalues in the Brillouin zone is shown in Figure~\ref{fig:band_structure2}.

\section{Discussions}\label{CFW}
We propose a finite element contour integral method for computing the resonances of metallic grating structures.
The finite element discretization allows for resolving the wave oscillation accurately when the
computational domain attains several length scales and the contrast of the medium coefficient is large. In addition, the multi-step contour integral technique solves the nonlinear eigenvalue problem successfully when the medium is dispersive. Numerical examples demonstrate the effectiveness of the proposed method.

The work on the computation of resonances in three-dimensional subwavelength structures is more challenging and will be reported elsewhere in the future. Note that the resulting nonlinear algebraic system can be highly ill-posed. Although {\bf Algorithm 1} is robust in the sense that the computed eigenvalues are accurate and close to the machine precision, there is still room to improve the efficiency of the algorithm. This will be explored in details in future.


The error analysis of the proposed method is beyond the scope of this work. Note that the error arises from several sources, including the finite element discretization, the truncated DtN map, and the multi-step eigensolver.

\section*{Acknowledgement}
Y. Xi is partially supported by the National Natural Science Foundation of China with Grant No.11901295, Natural Science Foundation of Jiangsu Province under BK20190431.
J. Lin is partially supported by the NSF grant DMS-2011148, and J. Sun is partially supported by an NSF Grant 2109949 and a Simons Foundation Collaboration Grant 711922.

\bibliography{references}

\end{document}